\documentclass[aos,seceqn,nameyear,dvips]{arximspdf}
\usepackage{graphics}
\doi{10.1214/09-AOS777}
\volume{38}
\issue{4}
\pubyear{2010}
\firstpage{1953}
\lastpage{1977}

\makeatletter

\newtheorem{prop}{Proposition}[section]
\newtheorem{lemma}[prop]{Lemma}
\newtheorem{cor}[prop]{Corollary}

\newtheorem{theorem}{Theorem}[section]

\newproclaim{Remark}{Remark}

\makeatother

\begin{document}
\begin{frontmatter}

\title{Inconsistency of bootstrap: The Grenander estimator}
\runtitle{Inconsistency of bootstrap: The Grenander estimator}

\begin{aug}
\author[A]{\fnms{Bodhisattva} \snm{Sen}\thanksref{t1}\ead[label=e1]{bs2528@columbia.edu}\ead[label=u1,url]{http://www.stat.columbia.edu/\texttildelow bodhi}\corref{}},
\author[B]{\fnms{Moulinath} \snm{Banerjee}\thanksref{t2}\ead[label=e2]{moulib@umich.edu}\ead[label=u2,url]{http://www.stat.lsa.umich.edu/\texttildelow moulib}}\\ and
\author[B]{\fnms{Michael} \snm{Woodroofe}\thanksref{t3}\ead[label=e3]{michaelw@umich.edu}\ead[label=u3,url]{http://www.stat.lsa.umich.edu/\texttildelow michaelw}}
\runauthor{B. Sen, M. Banerjee and M. Woodroofe}
\affiliation{Columbia University, University of Michigan and University
of Michigan}
\address[A]{B. Sen\\
Department of Statistics\\
Columbia University\\
1255 Amsterdam Avenue\\
New York, New York 10027\\
USA\\
\printead{e1}\\
\printead{u1}} %adresu isvedimo komanda gale!
\address[B]{M. Banerjee\\
M. Woodroofe\\
Department of Statistics\\
University of Michigan\\
1085 South University\\
Ann Arbor, Michigan 48109-1107\\
USA\\
\printead{e2}\\
\phantom{E-mail: }\printead*{e3}\\
\printead{u2}\\
\phantom{URL: }\printead*{u3}}
\end{aug}

\thankstext{t1}{Supported by NSF Grant DMS-09-06597.}
\thankstext{t2}{Supported by NSF Grant DMS-07-05288.}
\thankstext{t3}{Supported by NSF Grant AST-05-07453.}

% HISTORY:
\received{\smonth{10} \syear{2007}}
\revised{\smonth{10} \syear{2009}}

% ABSTRACT
%
\begin{abstract}
In this paper, we investigate the (in)-consistency of different
bootstrap methods for constructing confidence intervals in the class of
estimators that converge at rate $n^{1/3}$. The Grenander
estimator, the nonparametric maximum likelihood estimator of an unknown
nonincreasing density function $f$ on $[0,\infty)$, is a prototypical
example. We focus on this example and explore different approaches to
constructing bootstrap confidence intervals for $f(t_0)$, where $t_0
\in(0,\infty)$ is an interior point. We find that the bootstrap
estimate, when generating bootstrap samples from the empirical
distribution function $\mathbb{F}_n$ or its least concave majorant
$\tilde F_n$, does not have any weak limit in probability. We provide a
set of sufficient conditions for the consistency of any bootstrap
method in this example and show that bootstrapping from a smoothed
version of $\tilde F_n$ leads to strongly consistent estimators. The
$m$ out of $n$ bootstrap method is also shown to be consistent while
generating samples from $\mathbb{F}_n$ and $\tilde{F}_n$.
\end{abstract}

% KEYWORDS
%

\end{frontmatter}

%s1 ###
\section{Introduction}\label{intro}

If $X_1,X_2,\ldots,X_n \stackrel{\mathrm{ind}}{\sim}
f$ are a sample from a nonincreasing density $f$ on $[0,\infty)$,
then the Grenander estimator, the nonparametric maximum likelihood
estimator (NPMLE) $\tilde f_n$ of $f$ [obtained by maximizing the
likelihood $\prod_{i=1}^n f(X_i)$ over all nonincreasing densities],
may be described as follows: let $\mathbb{F}_n$ denote the empirical
distribution function (EDF) of the data, and $\tilde{F}_n$ its least
concave majorant. Then the NPMLE $\tilde{f}_n$ is the left-hand
derivative of $\tilde{F}_n$. This result is due to Grenander (\citeyear
{grenander56}) and
is described in detail by Robertson, Wright and Dykstra (\citeyear
{RWD88}), pages
326--328. Prakasa Rao (\citeyear{prakasa69}) obtained the asymptotic
distribution of
$\tilde{f}_n$, properly normalized: let $\mathbb{W}$ be a two-sided
standard Brownian motion on $\mathbb{R}$ with $\mathbb{W}(0) = 0$ and
\[
{\mathbb C} = \mathop{\arg\max}_{s\in{\mathbb R}} [{\mathbb W}(s) - s^2].
\]
If $0 < t_0 < \infty$ and $f'(t_0) \ne0$, then
%
%e1.1 ###
%
\begin{equation}
\label{eq:chrnff}
n^{1/3} \{ \tilde f_n(t_0) - f(t_0) \} \Rightarrow
2 \bigl|\tfrac{1}{2} f(t_0) f'(t_0) \bigr|^{1/3} {\mathbb C},
\end{equation}
where $\Rightarrow$ denotes convergence in distribution. There are
other estimators that exhibit similar asymptotic
properties; for example, Chernoff's (\citeyear{chernoff64}) estimator
of the mode, the
monotone regression estimator [Brunk (\citeyear{brunk70})],
Rousseeuw's (\citeyear{rousseeuw84}) least median of squares estimator,
and the estimator
of the shorth [Andrews et al. (\citeyear{andrewsetal72}) and Shorack
and Wellner (\citeyear{shorackWe86})].
The seminal paper by Kim and Pollard (\citeyear{kimPo90}) unifies
$n^{1/3}$-rate
of convergence problems in the more general \mbox{$M$-estimation} framework.
Tables and a survey of statistical problems in which the distribution
of ${\mathbb C}$ arises are provided by Groeneboom and Wellner
(\citeyear{GW01}).

The presence of nuisance parameters in the limiting distribution (\ref
{eq:chrnff}) complicates the construction of
confidence intervals. Bootstrap intervals avoid the problem of
estimating nuisance parameters and are generally reliable in problems
with $\sqrt{n}$ convergence rates. See Bickel and Freedman (\citeyear{BF81}),
Singh (\citeyear{singh81}), Shao and Tu (\citeyear{shaoTu95}) and its
references. Our aim in this
paper is to study the consistency of bootstrap methods for the
Grenander estimator with the hope that the monotone density estimation
problem will shed light on the behavior of bootstrap methods in similar
cube-root convergence problems.

There has been considerable recent interest in this question. Kosorok
(\citeyear{kosorok07}) show that bootstrapping from the EDF ${\mathbb
F}_n$ does not
lead to a consistent estimator of the distribution of $n^{1/3} \{
\tilde f_n(t_0) - f(t_0)\}$. Lee and Pun (\citeyear{leePun06}) explore
$m$ out of $n$
bootstrapping from the empirical distribution function in similar
nonstandard problems and prove the consistency of the method. L\'eger
and MacGibbon (\citeyear{legerMa06}) describe conditions for a
resampling procedure to
be consistent under cube root asymptotics and assert that these
conditions are generally not met while bootstrapping from the EDF. They
also propose a smoothed version of the bootstrap and show its
consistency for Chernoff's estimator of the mode. Abrevaya and Huang
(\citeyear{AH05}) show that bootstrapping from the EDF leads to inconsistent
estimators in the setup of Kim and Pollard (\citeyear{kimPo90}) and propose
corrections. Politis, Romano and Wolf (\citeyear{PRW99}) show that subsampling
based confidence intervals are consistent in this scenario.

Our work goes beyond that cited above as follows: we show that
bootstrapping from the NPMLE $\tilde{F}_n$ also leads to inconsistent
estimators, a result that we found more surprising, since $\tilde{F}_n$
has a density. Moreover, we find that \textit{the bootstrap estimator,
constructed from either the EDF or NPMLE, has no limit in probability}.
The finding is less than a mathematical proof, because one step in the
argument relies on simulation; but the simulations make our point
clearly. As described in Section~\ref{discussion}, our findings are
inconsistent with some claims of Abrevaya and Huang (\citeyear{AH05}).
Also, our
way of tackling the main issues differs from that of the existing
literature: we consider conditional distributions in more detail than
Kosorok (\citeyear{kosorok07}), who deduced inconsistency from
properties of
unconditional distributions; we directly appeal to the characterization
of the estimators and use a continuous mapping principle to deduce the
limiting distributions instead of using the ``switching'' argument [see
Groeneboom (\citeyear{groene85})] employed by Kosorok (\citeyear
{kosorok07}) and Abrevaya and Huang
(\citeyear{AH05}); and at a more technical level, we use the Hungarian
Representation theorem whereas most of the other authors use empirical
process techniques similar to those described by van der Vaart and
Wellner (\citeyear{VW00}).

Section \ref{prelim} contains a uniform version of (\ref{eq:chrnff})
that is used later on to study the consistency of different bootstrap
methods and may be of independent interest. The main results on
inconsistency are presented in Section \ref{boots_prob}. Sufficient
conditions for the consistency of a bootstrap method are presented in
Section \ref{smooth_boots} and applied to show that bootstrapping from
smoothed versions of $\tilde{F}_n$ does produce consistent estimators.
The $m$ out of $n$ bootstrapping procedure is investigated, when
generating bootstrap samples
from $\mathbb{F}_n$ and $\tilde{F}_n$. It is shown that both the
methods lead to consistent estimators under mild conditions on $m$. In
Section \ref{discussion}, we discuss our findings, especially the
nonconvergence and its implications. The
\hyperref[app]{Appendix}, provides the details of some arguments used
in proving the
main results.

%s2 ###
\section{Uniform convergence}\label{prelim}

For the rest of the paper, $F$ denotes a distribution function with
$F(0) = 0$ and a density $f$ that is nonincreasing on $[0,\infty)$ and
continuously differentiable near $t_0 \in(0,\infty)$ with nonzero
derivative $f'(t_0) < 0$. If $g\dvtx I \to{\mathbb R}$ is a bounded
function, write $\Vert g\Vert:= {\sup_{x\in I}} |g(x)|$. Next, let $F_n$
be distribution functions with $F_n(0) = 0$, that converge weakly to
$F$ and, therefore,
%
%e2.1 ###
%
\begin{equation}
\label{eq:cndtn0}
{\lim_{n\to\infty}} \| F_n-F \| = 0.
\end{equation}
Let $X_{n,1},X_{n,2},\ldots,X_{n,m_n} \stackrel{\mathrm{ind}}{\sim}
F_n$, where $m_n
\le n$ is a nondecreasing sequence of integers for which $m_n \to\infty
$; let ${\mathbb F}_{n,m_n}$ denote the EDF of $X_{n,1},X_{n,2},
\ldots,\break X_{n,m_n}$; and let
\[
\Delta_n := m_n^{1/3} \{\tilde f_{n,m_n}(t_0) - f_n(t_0)\},
\]
where $\tilde f_{n,m_n}(t_0)$ is the Grenander estimator computed from
$X_{n,1},X_{n,2},\ldots,\break X_{n,m_n}$ and $f_n(t_0)$ is the density of
$F_n$ at $t_0$ or a surrogate. Next, let $I_m = [-t_0m^{1/3},\infty
)$ and
%
%e2.2 ###
%
\begin{equation}\label{eq:Z_process}
\quad \mathbb{Z}_n(h) := m_n^{2/3} \{ \mathbb{F}_{n,m_n}(t_0 +
m_n^{-{1/3}}h) - \mathbb{F}_{n,m_n}(t_0) - f_n(t_0) m_n^{-{1/3}}h \}
\end{equation}
for $h \in I_{m_n}$ and observe that $\Delta_n$ is the left-hand
derivative at $0$ of the least concave majorant of ${\mathbb Z}_n$. It
is fairly easy to obtain the asymptotic distribution of $\mathbb{Z}_n$.
The asymptotic distribution of $\Delta_n$ may then be obtained from the
Continuous Mapping theorem. Stochastic processes are regarded as random
elements in $D(\mathbb{R})$, the space of right continuous functions on
$\mathbb{R}$ with left limits, equipped with the projection $\sigma
$-field and the topology of uniform convergence on compacta. See
Pollard (\citeyear{polllard84}), Chapters IV and V for background.

%s2.1 ###
\subsection{Convergence of ${\mathbb Z}_n$}

It is convenient to
decompose $\mathbb{Z}_n$ into the sum of $\mathbb{Z}_{n,1}$ and $\mathbb
{Z}_{n,2}$ where
\begin{eqnarray*}
\mathbb{Z}_{n,1}(h) &:=& m_n^{2/3} \{ (\mathbb
{F}_{n,m_n} - F_{n})(t_0 + m_n^{-{1/3}}h) - (\mathbb{F}_{n,m_n}-
F_{n})(t_0) \},\\
\mathbb{Z}_{n,2}(h) &:=& m_n^{2/3} \{
F_{n}(t_0 + m_n^{-{1/3}} h)- F_{n}(t_0) - f_n(t_0) m_n^{-{1/3}}h \}.
\end{eqnarray*}
Observe that ${\mathbb Z}_{n,2}$ depends only on $F_n$ and $f_n$; only
${\mathbb Z}_{n,1}$ depends on $X_{n,1},\ldots,\break X_{n,m_n}$. Let
$\mathbb{W}_1$ be a standard two-sided Brownian motion on $\mathbb{R}$
with $\mathbb{W}_1(0) = 0$, and $\mathbb{Z}_{1}(h) = \mathbb{W}_1 [f(t_0)h]$.
\begin{prop}\label{prop:Z_conv}
If
%
%e2.3 ###
%
\begin{equation}
\label{eq:cndtn1}
\lim_{n\to\infty} m_n^{1/3} | F_n(t_0+m_n^{-{1/3}}h) -
F_n(t_0) - f(t_0) m_n^{-{1/3}}h | = 0
\end{equation}
uniformly on compacts (in $h$), then ${\mathbb Z}_{n,1} \Rightarrow
{\mathbb Z}_1$; and if there is a continuous
function ${\mathbb Z}_2$ for which
%
%e2.4 ###
%
\begin{equation}
\label{eq:cndtn2}
\lim_{n\to\infty} {\mathbb Z}_{n,2}(h) = {\mathbb Z}_2(h)
\end{equation}
uniformly on compact intervals, then ${\mathbb Z}_n \Rightarrow
{\mathbb Z} := {\mathbb Z}_1 + {\mathbb Z}_2$.
\end{prop}
\begin{pf}
The Hungarian Embedding theorem of K\'omlos, Major and\break  Tusn\'ady
(\citeyear{kmt75})
is used. We may suppose that
$X_{n,i} = F_n^{\#}(U_{i})$, where $F_n^{\#}(u) = \inf\{x\dvtx F_n(x)
\ge
u\}$ and $U_{1},U_2,\ldots$ are i.i.d. Uniform$(0,1)$ random variables.
Let $\mathbb{U}_n$ denote the EDF of $U_{1},\ldots, U_{n}$, $\mathbb
{E}_n(t) = \sqrt{n} [\mathbb{U}_n(t) - t]$, and $\mathbb{V}_n = \sqrt
{m_n} (\mathbb{F}_{n,m_n} - F_n)$. Then $\mathbb{V}_n = \mathbb
{E}_{m_n} \circ F_n$. By Hungarian Embedding, we may also suppose that
the probability space supports a sequence of Brownian Bridges $\{\mathbb
{B}_n^0\}_{n \ge1}$ for which
%
%e2.5 ###
%
\begin{equation}
\label{eq:kmt}
{\sup_{0 \le t \le1}} | \mathbb{E}_{n}(t) - \mathbb{B}_n^0(t)| = O
\biggl[{\log(n)\over\sqrt{n}} \biggr] \qquad\mbox{a.s.},
\end{equation}
and a standard normal random variable $\eta$ that is independent of $\{
\mathbb{B}_n^0\}_{n \ge1}$. Define a version $\mathbb{B}_n$ of
Brownian motion by $\mathbb{B}_n(t) = \mathbb{B}_n^0(t) + \eta t$, for
$t \in[0,1]$. Then
%
%e2.6 ###
%
\begin{eqnarray}\label{eq:boots_proc1}
\mathbb{Z}_{n,1} (h) & = & m_n^{1/6} \{\mathbb
{E}_{m_n}[F_n(t_0 + m_n^{-{1/3}}h)] -
\mathbb{E}_{m_n}[F_{n}(t_0)] \} \nonumber\\[-8pt]\\[-8pt]
& = & m_n^{1/6} \{\mathbb{B}_{m_n}[F_n(t_0 + m_n^{-{1/3}}h)] -
\mathbb{B}_{m_n}[F_n(t_0)] \} + \mathbb{R}_{n}(h),\nonumber
\end{eqnarray}
where
\begin{eqnarray*}
|\mathbb{R}_n(h)| & \le& 2 m_n^{1/6} {\sup_{0
\le t \le1}} |\mathbb{E}_{m_n}(t) - \mathbb{B}_{m_n}^0(t)| \\
&&{} + m_n^{1/6}|\eta| |F_n(t_0 + m_n^{-{1/3}}h) - F_n(t_0)|
\rightarrow0
\end{eqnarray*}
uniformly on compacta w.p. 1 using (\ref{eq:cndtn1}) and (\ref
{eq:kmt}). Let
\[
\mathbb{X}_n(h) := m_n^{1/6}\{\mathbb
{B}_{m_n}[F_n(t_0 + m_n^{-{1/3}}h)] -
\mathbb{B}_{m_n}[F_n(t_0)] \}
\]
and observe that $\mathbb{X}_n$ is a mean zero Gaussian process defined
on $I_{m_n}$ with independent increments and
covariance kernel
\[
K_n(h_1,h_2) = m_n^{1/3} | F_n[t_0 + \operatorname{sign} \{h_1\}
m_n^{-{1/3}}(|h_1| \wedge|h_2|)] - F_n(t_0) | \mathbf{1}\{
h_1 h_2 > 0\}.
\]
It now follows from Theorem V.19 in Pollard (\citeyear{polllard84}) and
(\ref
{eq:cndtn1}) that $\mathbb{X}_n(h)$ converges in
distribution to $\mathbb{W}_1[f(t_0)h]$ in $D([-c,c])$ for every $c >
0$, establishing the first assertion of the
proposition. The second then follows from Slutsky's theorem.
\end{pf}

%s2.2 ###
\subsection{Convergence of $\Delta_n$}

Unfortunately, $\Delta_n$ is not
quite a continuous functional of ${\mathbb
Z}_n$. If $f\dvtx I \to{\mathbb R}$, write $f|J$ to denote the restriction
of $f$ to $J \subseteq I$; and if $I$ and $J$ are intervals and $f$ is
bounded, write $L_Jf$ for the least concave majorant of the
restriction. Thus, $\tilde{F}_n = L_{[0,\infty)}{\mathbb F}_n$ in the
\hyperref[intro]{Introduction}.
\begin{lemma}\label{lem:ww} Let $I$ be a closed interval; let $f\dvtx I
\to
{\mathbb R}$ be a bounded upper semi-continuous function on $I$; and
let $a_1,a_2,b_1,b_2 \in I$ with $b_1 < a_1 < a_2 < b_2$. If
$2f[{1\over2}(a_i+b_i)] > L_{I}
f(a_i) + L_{I} f(b_i), i = 1,2$, then $L_If(x) = L_{[b_1,b_2]}f(x)$ for
$a_1 \le x \le a_2$.
\end{lemma}
\begin{pf} This follows from the proof of Lemmas 5.1 and 5.2 of
Wang and Woodroofe (\citeyear{wangWoo07}). In that lemma
continuity was assumed, but only upper semi-continuity was used in the
(short) proof.
\end{pf}

Recall Marshall's lemma: if $I$ is an interval, $f: I \to{\mathbb R}$
is bounded, and $g\dvtx I \to{\mathbb R}$ is concave, then $\Vert
L_If-g\Vert\le\Vert f-g\Vert$. See, for example, Robertson, Wright
and Dykstra [(\citeyear{RWD88}), page 329] for a proof. Write $\tilde
{F}_{n,m_n} =
L_{[0,\infty)}{\mathbb F}_{n,m_n}$.
\begin{lemma}\label{lem:loc}
If $\delta> 0$ is so small that $F$ is strictly concave on
$[t_0-2\delta,t_0+2\delta]$ and (\ref{eq:cndtn0}) holds then ${\tilde
F}_{n,m_n} = L_{[t_0-2\delta,t_0+2\delta]}{\mathbb F}_{n,m_n}$ on
$[t_0-\delta,t_0+\delta]$ for all large $n$ w.p. 1.
\end{lemma}
\begin{pf} Since $F$ is strictly concave on $[t_0-2\delta,t_0+2\delta
], 2F(t_0\pm{3\over2}\delta) >
F(t_0\pm\delta) + F(t_0\pm2\delta)$. Then
\begin{eqnarray*}
\| \tilde{F}_{n,m_n} - F \| & \le& \|{\mathbb F}_{n,m_n} - F \|
\\
& \le& \|{\mathbb F}_{n,m_n} - F_n \| + \|F_n - F\| \\
& \le& \frac{1}{\sqrt{m_n}} \|\mathbb{E}_{m_n}\| + \|F_n - F \| \to0
\qquad\mbox{w.p. } 1
\end{eqnarray*}
by Marshall's lemma, (\ref{eq:cndtn0}) and the Glivenko--Cantelli
theorem. Thus,\break $2{\mathbb F}_{n,m_n}(t_0\pm
{3\over2}\delta) > \tilde{F}_{n,m_n}(t_0\pm\delta) + \tilde
{F}_{n,m_n}(t_0\pm2\delta)$, for all large $n$ w.p. 1, and Lemma \ref
{lem:loc} follows from Lemma \ref{lem:ww}.
\end{pf}
\begin{prop}\label{prop:loc} \textup{(i)} Suppose that (\ref{eq:cndtn0}) and
(\ref{eq:cndtn1}) hold and given $\gamma> 0$, there are $0 < \delta<
1$ and $C > 0$ for which
%
%e2.7 ###
%
\begin{equation}
\label{eq:cndtn3}
\bigl| F_n(t_0+h) - F_n(t_0) - f_n(t_0)h - \tfrac{1}{2} f'(t_0) h^2
\bigr| \le\gamma h^2 + C m_n^{-{2/3}}
\end{equation}
and
%
%e2.8 ###
%
\begin{equation}
\label{eq:cndtn4}
| F_n(t_0+h) - F_n(t_0) | \le C (|h| + m_n^{-{1/3}})
\end{equation}
for $|h| \le\delta$ and for all large $n$. If $J$ is a compact
interval and $\varepsilon>0$, then there is a compact $K \supseteq J$,
depending only on $\varepsilon, J, C, \gamma$, and $\delta$, for which
%
%e2.9 ###
%
\begin{equation}
\label{eq:loc}
P [L_{I_{m_n}}{\mathbb Z}_n = L_K{\mathbb Z}_n \mbox{ on } J ]
\ge1 - \varepsilon
\end{equation}
for all large $n$.

\textup{(ii)} Let $\mathbb{Y}$ be an a.s. continuous stochastic process on
$\mathbb{R}$ that is a.s. bounded above. If $\lim_{|h| \rightarrow
\infty} \mathbb{Y}(h)/|h| = -\infty$ a.e., then the compact $K
\supseteq J$ can be chosen so that
%
%e2.10 ###
%
\begin{equation}
\label{eq:loc2}
P [L_{\mathbb R}{\mathbb Y} = L_K{\mathbb Y} \mbox{ on } J ]
\ge1 - \varepsilon.
\end{equation}
\end{prop}
\begin{pf} For a fixed sequence ($F_n \equiv F$) (\ref{eq:loc}) would
follow from the assertion in Example 6.5 of Kim and Pollard (\citeyear
{kimPo90}), and
it is possible to adapt their argument to a triangular array using (\ref
{eq:cndtn3}) and
(\ref{eq:cndtn4}) in place of Taylor series expansion. A different
proof is presented in the \hyperref[app]{Appendix}.
\end{pf}

We will use the following easily verified fact. In its statement, the
metric space $\mathcal{X}$ is to be endowed with
the projection $\sigma$-field. See Pollard (\citeyear{polllard84}),
page 70.
\begin{lemma}\label{lemma:con_derv} Let $\{X_{n,c}\}, \{Y_n\}, \{W_c\}$
and $Y$ be sets of random elements taking values in a metric space
$(\mathcal{X}$,$d)$, $n=0,1,\ldots,$ and $c \in\mathbb{R}$. If for any
$\delta> 0$,
\begin{longlist}
\item $\lim_{c \rightarrow\infty} \limsup_{n \rightarrow\infty}
P\{d(X_{n,c},Y_n) > \delta\} = 0$,

\item $\lim_{c \rightarrow\infty} P\{d(W_{c},Y) > \delta\} = 0$,

\item $X_{n,c} \Rightarrow W_c$ as $n \rightarrow\infty$ for
every $c \in\mathbb{R}$,
\end{longlist}
then $Y_n \Rightarrow Y$ as $n \rightarrow\infty$.
\end{lemma}
\begin{cor}\label{cor:convdelta} If (\ref{eq:loc}) and (\ref{eq:loc2})
hold, and $\mathbb{Z}_n \Rightarrow\mathbb{Y}$, then
$L_{I_{m_n}}{\mathbb Z}_n \Rightarrow L_{\mathbb R}{\mathbb Y}$ in
$D({\mathbb R})$ and $ \Delta_n \Rightarrow(L_{\mathbb R}{\mathbb Y})'(0)$.
\end{cor}
\begin{pf} It suffices to show that $L_{I_{m_n}}{\mathbb Z}_n|J
\Rightarrow L_{\mathbb R}{\mathbb Y}|J$ in $D(J)$, for every compact
interval $J \subseteq\mathbb{R}$. Given $J$ and $\varepsilon> 0$, there
exists $K_\varepsilon$, a compact, $K_\varepsilon\supseteq J$, such that
(\ref{eq:loc}) and (\ref{eq:loc2}) hold. This verifies (i) and (ii)
of Lemma \ref{lemma:con_derv} with $c = 1/\varepsilon$, $X_{n,c}=
L_{K_\varepsilon}{\mathbb Z}_n$, $Y_n = L_{I_{m_n}}{\mathbb Z}_n$, $W_c =
L_{K_\varepsilon}{\mathbb Y}$, $Y = L_{\mathbb R}{\mathbb Y}$ and $d(x,y)
= \sup_{t \in J} |x(t) - y(t)|$. Clearly, $L_{K_\varepsilon}{\mathbb
Z}_n|J \Rightarrow L_{K_\varepsilon}{\mathbb Y}|J$ in $D(J)$, by the
Continuous Mapping theorem, verifying condition (iii). Thus,
$L_{I_{m_n}}{\mathbb Z}_n \Rightarrow L_{\mathbb R}{\mathbb Y}$ in
$D(\mathbb R)$. Another application of the Continuous Mapping theorem
[via the lemma on page 330 of Robertson, Wright and Dykstra (\citeyear
{RWD88})] in
conjunction with (\ref{eq:loc}), (\ref{eq:loc2}) and Lemma \ref
{lemma:con_derv} then shows that $\Delta_n = (L_{I_{m_n}}{\mathbb
Z}_n)'(0) \Rightarrow(L_{\mathbb R}{\mathbb Y})'(0)$.
\end{pf}
\begin{cor}\label{cor:convdelta2} If (\ref{eq:cndtn0}), (\ref
{eq:cndtn1}), (\ref{eq:cndtn2}), (\ref{eq:cndtn3}) and
(\ref{eq:cndtn4}) hold and
\[
\lim_{|h| \rightarrow\infty}
\mathbb{Z}(h)/|h| = -\infty,
\]
then
$L_{I_{m_n}}{\mathbb Z}_n \Rightarrow L_{\mathbb R}{\mathbb Z}$ in
$D({\mathbb R})$ and $ \Delta_n \Rightarrow
(L_{\mathbb R}{\mathbb Z})'(0)$; and if $\mathbb{Z}_{2}(h) = f'(t_0)
h^2/2$, then $ \Delta_n \Rightarrow2 |\frac{1}{2} f(t_0) f'(t_0)|^{1/3}
\mathbb C$, where $\mathbb C$ has Chernoff's distribution.
\end{cor}
\begin{pf} The convergence follows directly from Proposition \ref
{prop:loc} and\break Corollary~\ref{cor:convdelta}. Note that if $\mathbb
{Z}_{2}(h) = f'(t_0) h^2/2$, then (\ref{eq:loc}) and (\ref{eq:loc2})
hold\break and~Corollary~\ref{cor:convdelta} can be applied. That
$(L_{\mathbb R}{\mathbb Z})'(0)$ is distributed as\break $2 |\frac{1}{2}
f(t_0) f'(t_0)|^{1/3} \mathbb C$ when $\mathbb{Z}_{2}(h) =
f'(t_0) h^2/2$ follows from elementary properties of Brownian motion
via the ``switching'' argument of Groeneboom (\citeyear{groene85}).
\end{pf}

%s2.3 ###
\subsection{Remarks on the conditions}\label{remarks}

If $F_n \equiv F$
and $f_n \equiv f$, then clearly (\ref{eq:cndtn0}), (\ref{eq:cndtn1}),
(\ref{eq:cndtn2}), (\ref{eq:cndtn3}) and (\ref{eq:cndtn4}) all hold
with ${\mathbb Z}_2(h) = f'(t_0)h^2/2$ for some $0 < \delta< 1$ and $C
\ge f(t_0 - \delta)$ by a Taylor expansion of $F$ and the continuity of
$f$ and $f'$ around $t_0$.
\begin{cor}\label{cor:smoothF} If there is a $\delta> 0$ for which
$F_n$ has a continuously differentiable density
$f_n$ on $[t_0 - \delta, t_0 + \delta]$, and
%
%e2.11 ###
%
\begin{equation}
\label{eq:cndsmoothF}\quad
\lim_{n \rightarrow\infty} \Bigl[ \|
F_n - F\| + \sup_{|t - t_0| < \delta} \bigl(
|f_n(t) - f(t)| + |f_n'(t) - f'(t)| \bigr) \Bigr] = 0,
\end{equation}
then (\ref{eq:cndtn0}), (\ref{eq:cndtn1}), (\ref{eq:cndtn2}), (\ref
{eq:cndtn3}) and (\ref{eq:cndtn4}) hold with
${\mathbb Z}_2(h) = f'(t_0)h^2/2$, and $\Delta_n \Rightarrow2 |\frac
{1}{2} f(t_0) f'(t_0)|^{1/3} \mathbb C$.
\end{cor}
\begin{pf} The result can be immediately derived from Taylor expansion
of $F_n$ and the continuity of $f$ and $f'$
around $t_0$. To illustrate the idea, we show that (\ref{eq:cndtn3})
holds. Let $\gamma> 0$ be given. Clearly,
%
%e2.12 ###
%
\begin{eqnarray}\label{sm_boot_t6}
&&\biggl| F_n(t_0 + h) - F_n(t_0) - f_n(t_0) h - \frac{1}{2} h^2 f'(t_0)
\biggr| \nonumber\\[-8pt]\\[-8pt]
&&\qquad\le{\frac{1}{2} h^2 \sup_{|s| \le|h|}} |f_n'(t_0 + s) -
f'(t_0)|.\nonumber
\end{eqnarray}
Let $\delta> 0$ be so small that $|f'(t) - f'(t_0)| \le\gamma$ for
$|t - t_0| < \delta$, and let $n_0$ be so large that ${\sup_{|t - t_0|
\le\delta} }|f_n'(t) - f'(t)| \le\gamma$ for $n \ge n_0$. Then the
last line in
(\ref{sm_boot_t6}) is at most $\gamma h^2$ for $|h| \le\delta$ and $n
\ge n_0$.
\end{pf}

Another useful remark, used below, is that if $\lim_{n \rightarrow
\infty} m_n^{1/3} \|F_{m_n} - F\| = 0$, then (\ref{eq:cndtn0}),
(\ref{eq:cndtn1}) and (\ref{eq:cndtn4}) hold.

In the next three sections, we apply Proposition \ref{prop:Z_conv} and
Corollary \ref{cor:convdelta} to bootstrap samples drawn from the EDF,
its LCM, and smoothed versions thereof. Thus, let $X_1,X_2,\ldots
\stackrel{\mathrm{ind}}{\sim} F$; let ${\mathbb F}_n$ be the EDF of
$X_1,\ldots,X_n$; and
let $\tilde{F}_n$ be its LCM. If $F_n = {\mathbb F}_n$, then
(\ref{eq:cndtn0}), (\ref{eq:cndtn1}) and (\ref{eq:cndtn4}) hold almost
surely by the above remark, since
%
%e2.13 ###
%
\begin{equation}
\label{eq:lil}
\Vert{\mathbb F}_n - F\Vert= O \Biggl[\sqrt{\log\log(n)\over n}
\Biggr]\qquad\mbox{a.s.}
\end{equation}
by the Law of the Iterated Logarithm for the EDF, which may be deduced
from Hungarian Embedding; and the same is true if $F_n = \tilde{F}_n$
since $\Vert\tilde{F}_n-F\Vert\le\Vert{\mathbb F}_n-F\Vert$, by
Marshall's lemma.

If $m_n = n$ and $f_n = \tilde{f}_n$, then (\ref{eq:cndtn2}) is not
satisfied almost surely or in probability by either ${\mathbb F}_n$ or
$\tilde{F}_n$. For either choice, (\ref{eq:cndtn3}) is satisfied in
probability if $f_n = f$.
\begin{prop}\label{prop:edflcm}
Suppose that $m_n = n$ and that $f_n = f$. If $F_n$ is either the EDF
${\mathbb F}_n$ or its LCM $\tilde{F}_n$, then
for any $\gamma, \varepsilon> 0$, there are $C > 0$ and $0 < \delta< 1$
for which (\ref{eq:cndtn3}) holds with
probability at least $1- \varepsilon$ for all large $n$.
\end{prop}

The proof is included in the \hyperref[app]{Appendix}.

%s3 ###
\section{Inconsistency and nonconvergence of the bootstrap}\label{boots_prob}

We begin with a brief discussion of the bootstrap.

%s3.1 ###
\subsection{Generalities}\label{generalities}

Now, suppose that $X_1,X_2,\ldots\stackrel{\mathrm{ind}}{\sim} F$ are
defined on
a probability space $(\Omega,\mathcal{A}, P)$. Write ${\mathbf X}_n =
(X_1,\ldots,X_n)$ and suppose that the distribution function, $H_n$
say, of the random variable $R_n(\mathbf{X}_n,F)$ is of interest. The
bootstrap methodology can be broken into three simple steps:

\begin{longlist}
\item Construct an estimator $\hat{F}_n$ of $F$ from ${\mathbf X}_n$;

\item let $X_1^{*},\ldots,X_{m_n}^{*} \stackrel{\mathrm{ind}}{\sim
}\hat{F}_n$ be
conditionally i.i.d. given ${\mathbf X}_n$;

\item then let ${\mathbf X}_n^{*} = (X_1^{*},\ldots,X_{m_n}^{*})$
and estimate $H_n$ by the conditional
distribution function of $R_n^{*} = R({\mathbf X}_n^{*},\hat{F}_n)$
given ${\mathbf X}_n$; that is
\[
H_{n}^*(x) = P^*\{R^*_n \le x\},
\]
where $P^*\{\cdot\}$ is the conditional probability given the data
$\mathbf{X}_n$, or equivalently, the entire sequence $\mathbf{X} =
(X_1,X_2,\ldots)$.
\end{longlist}
Choices of $\hat{F}_n$ considered below are the EDF ${\mathbb F}_n$,
its least concave majorant
$\tilde{F}_n$, and smoothed versions thereof.

Let $d$ denote the Levy metric or any other metric metrizing weak
convergence of distribution functions. We say that
$H_{n}^*$ is \textit{weakly}, \textit{respectively}, \textit{strongly},
\textit{consistent} if
$d(H_n,H_n^*)\stackrel{P}{\rightarrow} 0$, respectively, $d(H_n,H_n^{*})
\to0$ a.s. If $H_{n}$ has a weak limit $H$, then consistency requires
$H_{n}^*$ to converge weakly to $H$, in probability; and if $H$ is
continuous, consistency requires
\[
{\sup_{x \in\mathbb{R}}} |H_{n}^*(x) - H(x)| \stackrel{P}{\rightarrow} 0
\qquad\mbox{as } n \rightarrow\infty.
\]
There is also the apparent possibility that $H_n^{*}$ could converge to
a random limit; that is, that there is a
$G\dvtx\Omega\times\mathbb{R} \to[0,1]$ for which $G(\omega,\cdot)$ is a
distribution function for each $\omega\in
\Omega$, $G(\cdot,x)$ is measurable for each $x \in\mathbb{R}$, and
$d(G,H_n^{*}) \stackrel{P}{\rightarrow} 0$. This
possibility is only apparent, however, if $\hat{F}_n$ depends only on
the order statistics. For if $h$ is a bounded
continuous function on $\mathbb{R}$, then any limit in probability of
$\int_{\mathbb{R}} h(x) H_n^*(\omega;dx)$ must be invariant under
finite permutations of $X_1,X_2,\ldots$ up to equivalence, and thus,
must be almost surely constant by the Hewitt--Savage zero--one law
[Breiman (\citeyear{breiman68})]. Let $\bar G(x) = \int_{\Omega}
G(\omega;x) P(d\omega
)$. Then $\bar G$ is a distribution function and $\int_{\mathbb{R}}
h(x) G(\omega;dx) = \int_{\mathbb{R}} h(x) \bar G(dx)$ a.s. for each
bounded continuous~$h$, and therefore for any countable collection of
bounded continuous $h$. It follows that $G(\omega;x) =\bar G(x)$ a.e.
$\omega$ for all $x$ by letting $h$ approach indicator functions.

Now let
\[
\Delta_n = n^{1/3} \{\tilde f_n(t_0) - f(t_0) \}
\quad\mbox{and}\quad \Delta_n^* = m_n^{1/3}
\{\tilde f_{n,m_n}^*(t_0) - \hat f_n(t_0) \},
\]
where $\hat f_n(t_0)$ is an estimate of $f(t_0)$, for example, $\tilde
f_n(t_0)$, and $\tilde f_{n,m_n}^*(t_0)$ is the
Grenander estimator computed from the bootstrap sample $X_1^{*},\ldots
,X_{m_n}^{*}$. Then weak (strong) consistency of the bootstrap means
%
%e3.1 ###
%
\begin{equation}
\label{eq:boot_limit}
{\sup_{x \in\mathbb R}} |P^*[\Delta_n^* \le x] - P[\Delta_n \le x]|
\rightarrow0
\end{equation}
in probability (almost surely), since the limiting distribution (\ref
{eq:chrnff}) of $\Delta_n$ is continuous.

%s3.2 ###
\subsection{Bootstrapping from the NPMLE $\tilde{F}_n$}\label{bootsNPMLE}

Consider now the case in which $m_n = n$, $\hat F_n = \tilde{F}_n$, and
$\hat{f}_n(t_0) = \tilde{f}_n(t_0)$. Let
\[
\mathbb{Z}_n^*(h) := n^{2/3} \{ \mathbb{F}_n^*(t_0 +
n^{-{1/3}}h) -
\mathbb{F}_n^*(t_0) - \tilde f_n(t_0)n^{-{1/3}}h \}
\]
for $h \in I_n = [-n^{1/3}t_0, \infty)$, where $\mathbb{F}_{n}^*$
is the EDF of the bootstrap sample
$X_1^*,\ldots,X_n^* \sim\tilde F_n$. Then $\mathbb{Z}_n^* = \mathbb
{Z}_{n,1}^* + \mathbb{Z}_{n,2}$, where
\begin{eqnarray}
\mathbb{Z}_{n,1}^* (h) &=& n^{2/3} \{ (\mathbb{F}_n^*-\tilde
{F}_n)(t_0 + n^{-{1/3}}h) -
(\mathbb{F}_n^*-\tilde{F}_n)(t_0) \},
\\
\mathbb{Z}_{n,2}(h) &=& n^{2/3} \{\tilde{F}_n(t_0 + h
n^{-{1/3}}) - \tilde{F}_n(t_0) - \tilde f_n(t_0) n^{-{1/3}}h \}.
\end{eqnarray}
Further, let $\mathbb{W}_1$ and $\mathbb{W}_2$ be two independent
two-sided standard Brownian motions on $\mathbb{R}$
with $\mathbb{W}_1(0) = \mathbb{W}_2(0) = 0$,
\begin{eqnarray*}
\mathbb{Z}_{1}(h) &=& \mathbb{W}_1[f(t_0)h],\\
\mathbb{Z}_{2}^0(h) &=& \mathbb{W}_2[f(t_0)h] +
\tfrac{1}{2}f'(t_0)h^2,\\
\mathbb{Z}_{2}(h) &=& L_{\mathbb R} \mathbb{Z}_{2}^0(h) - L_{\mathbb R}
\mathbb{Z}_{2}^0(0) - (L_{\mathbb{R}}
\mathbb{Z}_2^0)'(0)h,\\
\mathbb{Z} &=& \mathbb{Z}_1 + \mathbb{Z}_{2}.
\end{eqnarray*}
Then $\Delta_n^*$ equals the left derivative at $h = 0$ of the LCM of
$\mathbb{Z}_n^*$. It is first shown that
$\mathbb{Z}_n^*$ converges in distribution to $\mathbb{Z}$ but the
conditional distributions of $\mathbb{Z}_n^*$ do not have a limit. The
following two lemmas are needed.
\begin{lemma}\label{lemma:independence} Let $W_n$ and $W_n^*$ be random
vectors in $\mathbb{R}^l$ and $\mathbb{R}^k$,
respectively; let $Q$ and $Q^*$ denote distributions on the Borel sets
of $\mathbb{R}^l$ and $\mathbb{R}^k$; and let
$\mathcal{F}_n$ be sigma-fields for which $W_n$ is $\mathcal
{F}_n$-measurable. If the distribution of $W_n$ converges to $Q$ and
the conditional distribution of $W_n^*$ given $\mathcal{F}_n$ converges
in probability to $Q^*$, then the joint distribution of $(W_n,W_n^*)$
converges to the product measure $Q \times Q^*$.
\end{lemma}
\begin{pf} The above lemma can be proved easily using characteristic
functions. Kosorok (\citeyear{kosorok07}) includes a detailed proof.
\end{pf}

The next lemma uses a special case of the Convergence of Types theorem
[Lo\`{e}ve (\citeyear{loeve63}), page 203]: let $V, W, V_n$ be random
variables and
$b_n$ be constants; if $V$ has a nondegenerate distribution, $V_n
\Rightarrow V$ as $n \to
\infty$, and $V_n + b_n \Rightarrow W$, then $b = \lim_{n \to\infty}
b_n$ exists and $W$ has the same distribution as $V+b$.
\begin{lemma}\label{lemma:emp_boot1} Let $\mathbf{X}_n^*$ be a
bootstrap sample generated from the data $\mathbf{X}_n$. Let $Y_n :=
\psi_n(\mathbf{X}_n)$ and $Z_n := \phi_n(\mathbf{X}_n,\mathbf{X}_n^*)$
where $\psi_n\dvtx\mathbb{R}^n \to\mathbb{R}$ and $\phi_n\dvtx\mathbb{R}^{2
n} \to\mathbb{R}$ are measurable functions; and let $K_n$ and $L_n$ be
the conditional distribution functions of $Y_n + Z_n$ and $Z_n$ given
${\mathbf X}_n$, respectively. If there are distribution functions $K$
and $L$ for which $L$ is nondegenerate, $d(K_n,K) \stackrel
{P}{\rightarrow} 0$ and $d(L_n,L) \stackrel{P}{\rightarrow} 0$ then
there is a random variable $Y$ for which $Y_n \stackrel{P}{\rightarrow} Y$.
\end{lemma}
\begin{pf} If $\{n_k\}$ is any subsequence, then there exists a further
subsequence $\{n_{k_l}\}$ for which $d(K_{n_{k_l}},K) \rightarrow 0$
a.s. and $d(L_{n_{k_l}},L) \rightarrow 0$ a.s. Then $Y := \lim_{l
\rightarrow\infty} Y_{n_{k_l}}$ exists a.s. by the Convergence of
Types theorem,
applied conditionally given ${\mathbf X} := (X_1,X_2,\ldots)$ with
$b_l = Y_{n_{k_l}}$. Note that $Y$ does not depend on the subsequence
$n_{k_l}$, since two such subsequences can be joined to form another
subsequence using which we can argue the uniqueness.
\end{pf}
\begin{theorem}\label{thm:bootmle}
\textup{(i)} The conditional distribution of $\mathbb{Z}_{n,1}^*$ given
$\mathbf{X} = (X_1,\break X_2,\ldots)$ converges a.s. to the distribution
of $\mathbb{Z}_1$.

{\smallskipamount=0pt
\begin{longlist}[(iii)]
\item[(ii)] The unconditional distribution of $\mathbb{Z}_{n,2}$
converges to that of $\mathbb{Z}_2$ and the
unconditional distributions of $(\mathbb{Z}_{n,1}^*,\mathbb{Z}_{n,2})$,
and $\mathbb{Z}_n^*$ converge to those of
$(\mathbb{Z}_{1}, \mathbb{Z}_{2})$ and~$\mathbb{Z}$.

\item[(iii)] The unconditional distribution of $\Delta_n^*$ converges
to that of $(L_{\mathbb{R}}\mathbb{Z})'(0)$, and (\ref{eq:boot_limit}) fails.

\item[(iv)] Conditional on $\mathbf{X}$, the distribution of $\mathbb
{Z}_n^*$ does not have a weak limit in
probability.

\item[(v)] If the conditional distribution function of $\Delta_n^*$
converges in probability, then $(L_{\mathbb{R}} {\mathbb{Z}})'(0)$ and
${\mathbb Z}_2$ must be independent.
\end{longlist}}
\end{theorem}
\begin{pf} (i) The conditional convergence of $\mathbb
{Z}_{n,1}^*$ follows from Proposition \ref{prop:Z_conv} with $m_n = n$,
$F_n = \tilde{F}_n$, $\mathbb{F}_{n,m_n} = \mathbb{F}_n^*$, applied
conditionally given ${\mathbf X}$. It is only necessary to show that
(\ref{eq:cndtn1}) holds a.s., and this follows from the Law of the
Iterated Logarithm for ${\mathbb F}_n$ and Marshall's lemma, as
explained in Section \ref{remarks}. The unconditional limiting
distribution of
$\mathbb{Z}_{n,1}^*$ must also be that of $\mathbb{Z}_1$.

(ii) Let
\[
{\mathbb Z}_{n,2}^0(h) = n^{2/3}[{\mathbb F}_n(t_0+n^{-{1/3}}h) -
{\mathbb F}_n(t_0) - f(t_0)n^{-{1/3}}h]
\]
and observe that
\[
{\mathbb Z}_{n,2}(h) = L_{I_n} {\mathbb Z}_{n,2}^0(h) -
[L_{I_n}{\mathbb Z}_{n,2}^0(0) + (L_{I_n}{\mathbb Z}_{n,2}^0)'(0)h
].
\]
The unconditional convergence of ${\mathbb Z}_{n,2}^0$ and
$L_{I_n}{\mathbb Z}_{n,2}^0$ follow from Corollary \ref{cor:convdelta2}
applied with $F_n \equiv F$, as explained in Section \ref{remarks}. The
convergence in distribution of ${\mathbb Z}_{n,2}$ now follows from the
Continuous Mapping theorem, using Lemma \ref{lemma:con_derv} and
arguments similar to those in the proof of Corollary \ref{cor:convdelta}.

It remains to show that $\mathbb{Z}_{n,1}^*$ and $\mathbb{Z}_{n,2}^0$
are asymptotically independent, for example, the joint
limit distribution of $\mathbb{Z}_{n,1}^*$ and $\mathbb{Z}_{n,2}^0$ is
the product of their marginal limit
distributions. For this, it suffices to show that $(\mathbb
{Z}_{n,1}^*(t_1),\ldots,\break\mathbb{Z}_{n,1}^*(t_k))$ and
$(\mathbb{Z}_{n,2}^0(s_1),\ldots, \mathbb{Z}_{n,2}^0(s_l))$ are
asymptotically independent, for all choices $-\infty< t_1 < \cdots<
t_k <\infty$ and $-\infty< s_1 < \cdots< s_l <\infty$. This is an
easy consequence of
Lemma \ref{lemma:independence} applied with $W_n^* = (\mathbb
{Z}_{n,1}^*(t_1), \ldots,\mathbb{Z}_{n,1}^*(t_k))$ and $W_n = (\mathbb
{Z}_{n,2}^0(s_1),\ldots, \mathbb{Z}_{n,2}^0(s_l))$, and $\mathcal{F}_n
= \sigma(X_1,X_2,\ldots,X_n)$.\vspace*{1pt}

(iii) We will appeal to Corollary \ref{cor:convdelta} to find the
unconditional distribution of $\Delta_n^*$. We already know that
$\mathbb{Z}_n^*$ converges in distribution to $\mathbb{Z}$. That (\ref
{eq:loc2}) holds for the limit $\mathbb{Z}$ can be directly verified
from the definition of the process. We only have to show that (\ref
{eq:loc}) holds unconditionally with $\mathbb{Z}_n = \mathbb{Z}_n^*$.

Let $\varepsilon> 0$ and $\gamma> 0$ be given. By Proposition \ref
{prop:edflcm}, there exists $\delta> 0$ and $C>0$ such that $P(A_n)
\ge1- \varepsilon$ for all $n > N_0$, where
\[
A_n:= \bigl\{ \bigl|\tilde F_n(t_0 + h) + \tilde F_n(t_0) - f(t_0)h -
\tfrac{1}{2} f'(t_0) h^2 \bigr| \le\gamma h^2 + C n^{-{2/3}},
\forall|h| \le\delta\bigr\}.
\]
We can also assume that $|F(t_0 + h) + F(t_0) - f(t_0)h - (1/2) f'(t_0)
h^2| \le\gamma h^2$ for $|h| \le\delta$. Let ${\mathbb Y}_n^*(h) =
n^{2/3}[{\mathbb F}_n^{*}(t_0+n^{-{1/3}}h) - {\mathbb
F}_n^{*}(t_0) -
f(t_0)n^{-{1/3}}h]$, so that ${\mathbb Z}_n^*(h) = {\mathbb
Y}_n^*(h) - \Delta_n h$ for all $h \in I_n$, and
\[
L_K{\mathbb Z}^*_n = L_K{\mathbb Y}^*_n - \Delta_n h
\]
for all $h \in K$ for any interval $K \subseteq I_n$.

Let $G_n = \tilde F_n \mathbf{1}_{A_n} + F \mathbf{1}_{A_n^c}$ and let
$P_{G_n}^\infty$ denote the
probability when generating the bootstrap samples from $G_n$. Then
$G_n$ satisfies (\ref{eq:cndtn0}), (\ref{eq:cndtn1}), (\ref{eq:cndtn3})
and (\ref{eq:cndtn4}) a.s. with $m_n = n$, $F_n = G_n$, $\mathbb
{F}_{n,m_n} = \mathbb{F}^*_n \mathbf{1}_{A_n} + \mathbb{F}_n \mathbf
{1}_{A_n^c}$ and $f_n = f$. Let $J$ be a compact interval. By
Proposition \ref{prop:loc}, applied conditionally, there exists a
compact interval $K$ (not depending on $\omega$, by the \textit{remark}
near the end of the proof of Proposition~\ref{prop:loc}) such that $K
\supseteq J$ and
\[
P_{G_n}^\infty[L_{I_n}{\mathbb Y}^*_n = L_{K} {\mathbb Y}^*_n \mbox{
on } J] (\omega) \ge1-\varepsilon
\]
for $n \ge N(\omega) $ for a.e. $\omega$. As $N(\cdot)$ is bounded in
probability, there exists $N_1 > 0$ such that $P(B) \ge1 - \varepsilon$,
where $B := \{\omega\dvtx N(\omega) \le N_1\}$. By increasing $N_1$ if
necessary, let us also suppose that $N_1 \ge N_0$. Then
\begin{eqnarray*}
P[L_{I_{m_n}} \mathbb{Z}_n^* = L_K \mathbb{Z}_n^* \mbox{ on } J] & = &
P[L_{I_{m_n}} \mathbb{Y}_n^* = L_K \mathbb{Y}_n^* \mbox{ on } J]
\\
& \ge& \int_{A_n} P^*[L_{I_{m_n}} \mathbb{Y}_n^* = L_K \mathbb{Y}_n^*
\mbox{ on } J](\omega) \,dP(\omega) \\
& = & \int_{A_n} P_{G_n}^\infty[L_{I_{m_n}} \mathbb{Y}_n^* = L_K \mathbb
{Y}_n^* \mbox{ on } J](\omega) \,dP(\omega) \\
& \ge& \int_{A_n \cap B} P_{G_n}^\infty[L_{I_{m_n}} \mathbb{Y}_n^* =
L_K \mathbb{Y}_n^* \mbox{ on } J](\omega) \,dP(\omega) \\
& \ge& \int_{A_n \cap B} (1 - \varepsilon) \,dP(\omega) \ge1 - 3
\varepsilon\qquad \mbox{for all } n \ge N_1
\end{eqnarray*}
as $P(A_n \cap B) \ge1 - 2 \varepsilon$ for $n \ge N_1$. Thus, (\ref
{eq:loc}) holds and Corollary \ref{cor:convdelta} gives $\Delta_n^*
\Rightarrow(L_{\mathbb{R}} \mathbb{Z})'(0)$.

If (\ref{eq:boot_limit}) holds in probability, then the unconditional
limit distribution of $\Delta_n^{*}$ would be that of $2 |\frac{1}{2}
f(t_0) f'(t_0)|^{1/3} \mathbb C$, which is different from the
distribution of $(L_{\mathbb R}{\mathbb Z})'(0)$, giving rise to a
contradiction.

(iv) We use the method of contradiction. Let $Z_n := \mathbb
{Z}_{n,1}^*(h_0)$ and $Y_n := \mathbb{Z}_{n,2}(h_0)$ for some fixed
$h_0 > 0$ (say $h_0 = 1$) and suppose that the conditional distribution
function of $Z_n + Y_n = \mathbb{Z}_n^*(h_0)$ converges in probability
to the distribution function~$G$. By Proposition \ref{prop:Z_conv}, the
conditional distribution of $Z_n$ converges in probability to a normal
distribution, which is obviously nondegenerate. Thus, the assumptions
of Lemma \ref{lemma:emp_boot1} are satisfied and we conclude that $Y_n
\stackrel{P}{\rightarrow} Y$, for some random variable $Y$. It then
follows from the Hewitt--Savage zero--one law that $Y$ is a constant, say
$Y = c_0$ w.p. 1. The contradiction arises since $Y_n$ converges in
distribution to $\mathbb{Z}_2(h_0)$ which is not a constant a.s.

%f1 ###
%
\begin{figure}

\includegraphics{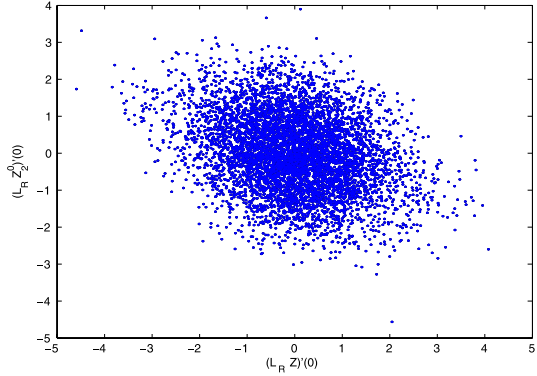}

\caption{Scatter plot of $10\mbox{,}000$ random draws of $((L_{\mathbb
R}{\mathbb Z})'(0),(L_{\mathbb R}{\mathbb
Z}_2^{0})'(0))$ when $f(t_0) = 1$ and $f'(t_0) = -2$.}\label{fig:simul}
\end{figure}

(v) We can show that the (unconditional) joint distribution of
$(\Delta_n^*, \mathbb{Z}_{n,2}^0)$ converges to that of
$((L_{\mathbb{R}} {\mathbb{Z}})'(0), {\mathbb Z}_2^0)$. But
$\Delta_n^*$ and $\mathbb{Z}_{n,2}^0$ are asymptotically independent by
Lemma \ref{lemma:independence} applied to $W_n =
(\mathbb{Z}_{n,2}^0(t_1),\mathbb{Z}_{n,2}^0(t_2),
\ldots,\mathbb{Z}_{n,2}^0(t_l))$, where $t_i \in\mathbb{R}$, $W_n^* =
\Delta_n^*$ and $\mathcal{F}_n = \sigma(X_1,X_2,\ldots,X_n)$. Therefore,
$(L_{\mathbb{R}} {\mathbb{Z}})'(0)$ and ${\mathbb Z}_2^0$ are
independent. The proposition follows directly since $\mathbb{Z}_2$ is a
measurable function of $\mathbb{Z}_{2}^0$.
%\rightqed
\end{pf}

If the conditional distribution of $\Delta^*_n$ converges in
probability, as a consequence of (v) of Theorem
\ref{thm:bootmle}, $(L_{\mathbb R}{\mathbb Z})'(0)$ and $(L_{\mathbb
R}{\mathbb Z}_2^{0})'(0)$ must also be independent. Figure \ref
{fig:simul} shows the scatter plot of $(L_{\mathbb R}{\mathbb Z})'(0)$
and $(L_{\mathbb R}{\mathbb Z}_2^{0})'(0)$ obtained from a simulation
study with $10\mbox{,}000$ samples, $f(t_0) = 1$ and $f'(t_0) = -2$. The
correlation coefficient obtained $-0.2999$ is highly significant
($p$-value $< 0.0001$). Thus, when combined with simulations, (v) of
Theorem \ref{thm:bootmle} strongly suggests that the conditional
distribution of $\Delta_n^*$ does not converge in probability.

%s3.3 ###
\subsection{Bootstrapping from the EDF}

A similar, slightly simpler
pattern arises if the bootstrap sample is drawn
from $\hat{F}_n = {\mathbb F}_n$. Define ${\mathbb Z}_n^*$ as before,
and let $\mathbb{Z}_{n,1}^* (h) = n^{2/3} \{
(\mathbb{F}_n^*-{\mathbb F}_n)(t_0 + n^{-{1/3}}h) - (\mathbb
{F}_n^*- {\mathbb F}_n)(t_0)\}$ and
$\mathbb{Z}_{n,2}(h) = n^{2/3} \{{\mathbb F}_n(t_0 + h n^{-{1/3}}) -
{\mathbb F}_n(t_0) - \tilde f_n(t_0)
n^{-{1/3}}h\}$. Then ${\mathbb Z}_n^* = {\mathbb Z}_{n,1}^* +
{\mathbb Z}_{n,2}$. Recall the definition of the
processes $\mathbb{W}_1$, $\mathbb{W}_2$, $\mathbb{Z}_1$, $\mathbb
{Z}^0_2$ in Section \ref{bootsNPMLE}. Define
\[
\mathbb{Z}_2(h) = \mathbb{Z}_2^0(h) -(L_{\mathbb{R}} \mathbb
{Z}_2^0)'(0) h.
\]
\begin{theorem}\label{thm:bootmle2}
\textup{(i)} The conditional distribution of $\mathbb{Z}_{n,1}^*$ given
$\mathbf{X} = (X_1,X_2,\break\ldots)$ converges a.s. to the distribution
of $\mathbb{Z}_1$.

{\smallskipamount=0pt
\begin{longlist}[(iii)]
\item[(ii)] The unconditional distribution of $\mathbb{Z}_{n,2}$
converges to that of $\mathbb{Z}_2$ and the
unconditional distributions of $(\mathbb{Z}_{n,1}^{*},\mathbb
{Z}_{n,2})$, and $\mathbb{Z}_n^*$ converge to those of $(\mathbb
{Z}_{1},\mathbb{Z}_{2})$ and~$\mathbb{Z}$.

\item[(iii)] The unconditional distribution of $\Delta_n^*$ converges
to that of $(L_{\mathbb{R}}\mathbb{Z})'(0)$, and (\ref{eq:boot_limit}) fails.

\item[(iv)] Conditional on $\mathbf{X}$, the distribution of $\mathbb
{Z}_n^*$ does not have a weak limit in
probability.

\item[(v)] If the conditional distribution function of $\Delta_n^*$
converges in probability, then
$(L_{\mathbb{R}}{\mathbb{Z}})'(0)$ and ${\mathbb Z}_2$ must be independent.
\end{longlist}}
\end{theorem}
\begin{Remark*} The proof of this theorem runs along similar lines to
that of Theorem \ref{thm:bootmle}. We briefly
highlight the differences.
\end{Remark*}

\begin{longlist}
\item The conditional convergence of $\mathbb{Z}_{n,1}^*$ follows
from Proposition \ref{prop:Z_conv} with $m_n = n$, $F_n = \mathbb
{F}_n$, $\mathbb{F}_{n,m_n} = \mathbb{F}_n^*$, applied conditionally
given ${\mathbf X}$. It is only necessary to show that (\ref
{eq:cndtn1}) is satisfied almost surely, and this follows from the Law
of the Iterated Logarithm for ${\mathbb F}_n$, as explained in
Section \ref{remarks}. Then the unconditional limiting distribution of
$\mathbb
{Z}_{n,1}^*$ must also be that of $\mathbb{Z}_1$.

\item The proof is similar to that of (ii) of Theorem \ref
{thm:bootmle}, except that now $\mathbb{Z}_{n,2}(h) = \mathbb
{Z}^0_{n,2}(h) - (L_{I_n} \mathbb{Z}^0_{n,2})'(0) h$.
\end{longlist}

The proofs of (iii)--(v) are very similar to that of (iii)--(v) of
Theorem \ref{thm:bootmle}.

%s3.4 ###
\subsection{Performance of the bootstrap methods in finite
samples}\label{simul}

In this subsection, we illustrate the poor finite sample performance
of the two inconsistent bootstrap schemes, namely, bootstrapping from
the EDF $\mathbb{F}_n$ and the NPMLE $\tilde F_n$. Table~\ref
{PerfBootsMeth} shows the estimated coverage probabilities of nominal
95\% confidence intervals for $f(1)$ using the two bootstrap methods
%
%t1 ###
%
\begin{table}[b]
\caption{Estimated coverage probabilities of nominal 95\% confidence
intervals for $f(1)$ while bootstrapping from the EDF $\mathbb{F}_n$
and NPMLE $\tilde F_n$, with varying sample size $n$ for the two
models: $\operatorname{Exponential}(1)$ (left) and $|Z|$ where $Z \sim
\operatorname{Normal}(0,1)$ (right)}\label{PerfBootsMeth}
\begin{tabular*}{\tablewidth}{@{\extracolsep{\fill}}lccccc@{}}
\hline
$\bolds n$ & \textbf{EDF} & \textbf{NPMLE} & $\bolds n$ & \textbf{EDF} & \textbf{NPMLE} \\
\hline
\phantom{0}$50$ & 0.747 & 0.720 & \phantom{0}$50$ & 0.761 & 0.739 \\
$100$ & 0.776 & 0.755  &  $100$ & 0.778 & 0.757 \\
$200$ & 0.802 & 0.780  &  $200$ & 0.780 & 0.762 \\
$500$ & 0.832 & 0.797  &  $500$ & 0.788 & 0.755 \\
\hline
\end{tabular*}
%
%$\bolds n$ & \textbf{EDF} & \textbf{NPMLE} \\
%$50$ & 0.761 & 0.739 \\
%$100$ & 0.778 & 0.757 \\
%$200$ & 0.780 & 0.762 \\
%$500$ & 0.788 & 0.755 \\
\end{table}
for different sample sizes, when the true distribution is assumed to be
Exponential(1) and $|\mathrm{Normal}(0,1)|$, respectively. We used 1000
bootstrap samples to compute each confidence interval and then
constructed 1000 such confidence intervals to estimate the actual
coverage probabilities. As is clear from the table the coverage
probabilities fall well short of the nominal 0.95 value. Leger and
MacGibbon (\citeyear{legerMa06}) also illustrate such a discrepancy in
the nominal and
actual coverage probabilities while bootstrapping from the EDF for the
Chernoff's estimator of the mode.

\begin{figure}

\includegraphics{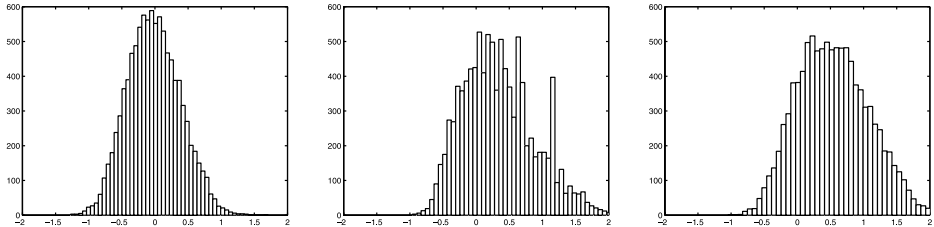}

\caption{Histograms of the exact distribution of $\Delta_n$ (left
panel) and the two bootstrap distributions while drawing bootstrap
samples from $\mathbb{F}_n$ (middle panel) and $\tilde F_n$ (right
panel) for $n = 500$.}\label{fig:InConsEDF_Hist}
\end{figure}

Figure \ref{fig:InConsEDF_Hist} shows the histograms (computed from
10,000 bootstrap samples) of the two inconsistent bootstrap
distributions obtained from a single sample of 500 Exponential(1)
random variables along with the histogram of the exact distribution of
$\Delta_n$ (obtained from simulation). The bootstrap distributions are
skewed and have very different shapes and supports compared to that on
the left panel of Figure \ref{fig:InConsEDF_Hist}. The histograms
illustrate the inconsistency of the bootstrap procedures.

%
%f2 ###
%

%
%
%f3 ###
%
\begin{figure}[b]

\includegraphics{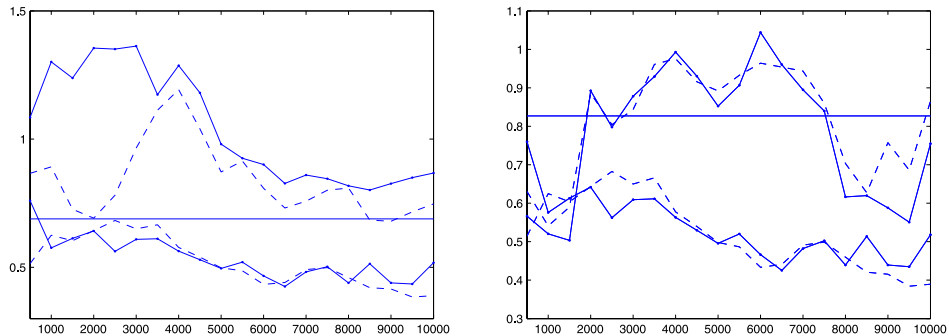}

\caption{Estimated 0.95 quantile of the bootstrap distribution while
generating the bootstrap samples from $\mathbb{F}_n$ (dashed lines) and
$\tilde F_n$ (solid-dotted lines) for two independent data sequences
along with the 0.95 quantile of the limit distribution of $\Delta_n$
(solid line) for the two models: $\operatorname{Exponential}(1)$ (left panel) and $|Z|$
where $Z \sim\operatorname{Normal}(0,1)$ (right panel).}\label{fig:InConsQuantile}
\end{figure}

The estimated coverage probabilities in Table \ref{PerfBootsMeth} are
unconditional [see (iii) of Theorems \ref{thm:bootmle} and \ref
{thm:bootmle2}] and do not provide direct evidence to suggest that the
conditional distribution of $\Delta_n^*$ does not converge in
probability. Figure \ref{fig:InConsQuantile} shows the estimated 0.95
quantile of the bootstrap distribution for two independent data
sequences as the sample size increases from 500 to 10,000, for the two
bootstrap procedures, and for both the models (exponential and normal).
The bootstrap quantile fluctuates enormously even at very large sample
sizes and shows signs of nonconvergence. If the bootstrap were
consistent, the estimated quantiles should converge to 0.6887 (0.8269),
the 0.95 quantile of the limit distribution of $\Delta_n$, indicated by
the solid line in Figure \ref{fig:InConsQuantile}. From the left panel
of Figure \ref{fig:InConsQuantile}, we see that the estimated bootstrap
0.95 quantiles (obtained from the two procedures) for one data sequence
stays below 0.6887, while for the other, the 0.95 quantiles stay above
0.6887, indicating the strong dependence on the sample path. Note that
if the bootstrap distributions had a limit, then Figure \ref
{fig:InConsQuantile} suggests that the limit varies with the sample
path, and that is impossible as explained in Section \ref
{generalities}. This provides evidence for the nonconvergence of the
bootstrap estimator.

%s4 ###
\section{Consistent bootstrap methods}\label{smooth_boots}

The main reason for the inconsistency of bootstrap methods discussed in
the previous section is the lack of smoothness of the distribution
function from which the bootstrap samples are generated. The EDF
$\mathbb{F}_n$ does not have a density, and $\tilde F_n$ does not have
a differentiable density, whereas $F$ is assumed to have a nonzero
differentiable density at $t_0$. At a more technical level, the lack of
smoothness manifests itself through the failure of (\ref{eq:cndtn2}).

The results from Section \ref{prelim} may be directly applied to derive
sufficient conditions on the smoothness of the distribution from which
the bootstrap samples are generated. Let
$X_1,X_2,\ldots\stackrel{\mathrm{ind}}{\sim}
F$; let $\hat F_n$ be an estimate of $F$ computed from $X_1,\ldots
,X_n$; and let $\hat f_n$ be the density of $\hat F_n$ or a surrogate,
as in Section \ref{boots_prob}.
\begin{theorem}\label{thm:cons_sm_boot} If (\ref{eq:cndtn0}), (\ref
{eq:cndtn1}), (\ref{eq:cndtn2}), (\ref{eq:cndtn3}) and
(\ref{eq:cndtn4}) hold a.s. with $F_n = \hat F_n$ and $f_n = \hat f_n$,
then the bootstrap estimate is strongly
consistent, for example, (\ref{eq:boot_limit}) holds w.p. 1. In
particular, the bootstrap estimate is strongly consistent if
there is a $\delta> 0$ for which $\hat F_n$ has a continuously
differentiable density $\hat f_n$ on $[t_0 - \delta,t_0 + \delta]$, and
(\ref{eq:cndsmoothF}) holds a.s. with $F_n = \hat F_n$ and $f_n = \hat f_n$.
\end{theorem}
\begin{pf} That $\Delta_n^*$ converges weakly to the distribution on
the right-hand side of (\ref{eq:chrnff}) a.s. follows from Corollary
\ref
{cor:convdelta2} applied conditionally given $\mathbf{X}$ with $F_n =
\hat F_n$ and $f_n = \hat f_n$. The second assertion follows similarly
from Corollary \ref{cor:smoothF}.
\end{pf}

%s4.1 ###
\subsection{Smoothing $\tilde F_n$}

We show that generating bootstrap
samples from a suitably smoothed version of
$\tilde{F}_n$ leads to a consistent bootstrap procedure. To avoid
boundary effects and ensure that the smoothed version has a decreasing
density on $(0,\infty)$, we use a logarithmic transformation. Let $K$
be a twice continuously differentiable symmetric density for which
%
%e4.1 ###
%
\begin{equation}
\label{eq:kay}
\int_{-\infty}^{\infty} [K(z)+|K'(z)|+|K''(z)|]e^{\eta|z|}\,dz < \infty
\end{equation}
for some $\eta> 0$. Let
%
%e4.2 ###
%
\begin{eqnarray}\label{eq:fchck1}
K_h(x,u) & = & {1\over hx}K \biggl[{1\over h}\log\biggl({u\over x}\biggr) \biggr]\quad
\mbox{and} \nonumber\\[-8pt]\\[-8pt]
\check{f}_n(x) & = & \int_{0}^{\infty} K_h(x,u)\tilde{f}_n(u)\,du = \int
_{0}^{\infty} K_h(1,u)\tilde{f}_n(xu)\,du.\nonumber
\end{eqnarray}
Thus, $e^y\check{f}_n(e^y) = \int_{-\infty}^{\infty}
h^{-1}K[h^{-1}(y-z)] \tilde{f}_n(e^z)e^z \,dz$. Integrating and using
capital letters to denote distribution functions,
\begin{eqnarray*}
\label{eq:fchck2}
\check{F}_n(e^y) & = & \int_{-\infty}^y \check f_n(e^s) e^s \,ds
\\
& = & \int_{-\infty}^y \int_{-\infty}^\infty\frac{1}{h} K \biggl( \frac
{s - v}{h} \biggr) \tilde f_n(e^v) e^v \,dv \,ds
\\
& = & \int_{-\infty}^{\infty} K(z)\tilde{F}_n(e^{y-hz})\,dz.
\end{eqnarray*}
Alternatively, integrating (\ref{eq:fchck1}) by parts yields
\[
\check{f}_n(x) = -\int_{0}^{\infty} {
\partial\over
\partial u}K_h(x,u)\tilde{F}_n(u)\,du.
\]
The proof of (\ref{eq:boot_limit}) requires showing that $\check{F}_n$
and its derivatives are sufficiently close to
those of $F$, and it is convenient to separate the estimation error
$\check{F}_n-F$ into sampling and approximation
error. Thus, let
%
%e4.3 ###
%
\begin{equation}
\label{eq:fbar}
\bar{F}_h(e^y) = \int_{-\infty}^{\infty} K(z) F(e^{y-hz})\,dz.
\end{equation}
We denote the first and second derivatives of $\bar{F}_h$ by $\bar
{f}_h$ and $\bar{f}_h'$, respectively. Recall that $F$ is assumed to
have a nonincreasing density on $(0,\infty)$ that is continuously
differentiable near $t_0$.
\begin{lemma}\label{lem:fbar} ${\lim_{h \rightarrow0}} \| \bar F_h - F\|
= 0$, and there is a $\delta> 0$ for which
%
%e4.4 ###
%
\begin{equation}
\label{eq:fbar2}
\lim_{h\to0} \sup_{|x-t_0|\le\delta} [ |\bar{f}_h(x)-f(x)| +
|\bar{f}_h'(x)-f'(x)| ] = 0.
\end{equation}
\end{lemma}
\begin{pf} First, observe that
\[
\bar{F}_h(e^y)-F(e^y) = \int_{-\infty}^{\infty}
K(z)[F(e^{y-hz})-F(e^y)]\,dz
\]
by (\ref{eq:fbar}). That $\lim_{h \rightarrow0} \bar F_h(x) = F(x)$
for all $x \ge0$ follows easily from the Dominated Convergence
theorem, and uniform convergence then follows from Polya's theorem.
This establishes the first assertion of the lemma. Next, consider (\ref
{eq:fbar2}). Given $t_0 > 0$, let $y_0 = \log(t_0)$ and let $\delta>
0$ be so small that $e^yf(e^y)$ is continuously differentiable (in $y$)
on $[y_0-2\delta,y_0+2\delta]$. Then
\begin{eqnarray*}
\bar{f}_h(x) - f(x) & = & \int_{-\infty}^{\infty} K(z)[f(xe^{hz})-f(x)]
e^{h z} \,dz \\
&&{} + f(x) \int_{-\infty}^{\infty} (e^{hz} - 1) K(z) \,dz
\end{eqnarray*}
and thus
\[
{\sup_{|x-t_0|\le\delta}} |\bar{f}_h(x)-f(x)| \le{\int_{-\infty}^{\infty
} \sup_{|x-t_0|\le\delta}} |f(xe^{hz}) - f(x)|
e^{h z} K(z)\,dz + O(h^2)
\]
for any $0 < \delta< t_0$. For sufficiently small $\delta$, the
integrand approach zero as $h \to0$; and it is bounded by $\sup
_{|x-t_0|\le\delta} (e^{-hz}/x + f(x)) e^{hz} K(z)$, since $f(x) \le
1/x$ for all $x > 0$. So the right-hand side approaches zero as $h \to
0$ by
the Dominated Convergence theorem. That ${\sup_{|x-t_0|\le\delta}} |\bar
{f}_h'(x)-f'(x)| \to0$ may be established similarly.
\end{pf}
\begin{theorem} Let $K$ be a twice continuously differentiable,
symmetric density for which (\ref{eq:kay}) holds. If
\[
h = h_n \to0 \quad\mbox{and}\quad h_n^{2}\sqrt{n\over\log\log(n)}
\to\infty,
\]
then the bootstrap estimator is strongly consistent; that is, (\ref
{eq:boot_limit}) holds a.s.
\end{theorem}
\begin{pf} By Theorem\vspace*{1pt} \ref{thm:cons_sm_boot}, it suffices to show that
(\ref{eq:cndsmoothF}) holds a.s. with $\hat
F_n = \check F_n$ and $\hat f_n = \check f_n$; and this would follow
from
\[
\|\check F_n - \bar F_h\| + \sup_{|x - t_0| \le\delta} [|\check f_n(x)
- \bar f_h(x)| + |\check f_n'(x) - \bar
f_h'(x)|] \rightarrow0 \qquad\mbox{a.s.}
\]
for some $\delta> 0$ and Lemma \ref{lem:fbar}. Clearly, using (\ref
{eq:fchck2}),
%e4.5 ###
%
\begin{equation}\label{eq:Fn-Fh}
\check{F}_n(e^y)-\bar{F}_h(e^y) = \frac{1}{h} \int_{-\infty}^{\infty}
[\tilde{F}_n(e^{t}) - F(e^{t})] K \biggl( \frac{y
- t}{h} \biggr) \,dt
\end{equation}
for all $y$, so that
\[
\| \check{F}_n - \bar F_h \| \le\|{\tilde F}_n - F \| \le\|{\mathbb
F}_n - F \| = O \bigl[ \sqrt{\log\log(n)/n}
\bigr] \qquad\mbox{a.s.}
\]
by Marshall's lemma and the Law of the Iterated Logarithm.
Differentiating (\ref{eq:Fn-Fh}) gives
\[
\check f_n(e^{y}) - \bar f_h(e^{y}) = \frac{e^{-y}}{h^2} \int_{-\infty
}^{\infty} [\tilde{F}_n(e^{t}) - F(e^{t})]
K' \biggl( \frac{y - t}{h} \biggr) \,dt.
\]
Differentiating (\ref{eq:Fn-Fh}) again and then taking absolute values
and considering $0 < h \le1$, we get
\begin{eqnarray*}
&&\hspace*{-6pt} \sup_{|x - t_0| \le\delta} \{|\check f_n(x) - \bar f_h(x)| +
|\check f_n'(x) - \bar f_h'(x)| \} \\
&&\hspace*{-6pt}\qquad \le \frac{M}{h^3} \sup_{|x - t_0| \le\delta} \int_{-\infty
}^{\infty} |\tilde{F}_n(e^{t}) - F(e^{t})| \biggl[ \biggl| K' \biggl( \frac
{\log x -
t}{h} \biggr) \biggr| + \biggl| K'' \biggl( \frac{\log x - t}{h} \biggr)
\biggr| \biggr] \,dt \\
&&\hspace*{-6pt}\qquad \le \frac{M}{h^2} \|{\mathbb F}_n - F \| \int_{-\infty}^{\infty}
[|K'(z)| + |K''(z)|] \,dz
\rightarrow 0 \qquad\mbox{a.s.}
\end{eqnarray*}
for a constant $M > 0$, as $h_n^2 \sqrt{n/ \log\log(n)} \rightarrow
\infty$, where Marshall's lemma and the Law of Iterated Logarithm have
been used again.
\end{pf}

%s4.2 ###
\subsection{$m$ out of $n$ bootstrap}\label{m_n_boots}

In Section \ref
{boots_prob}, we showed that the two most intuitive methods of
bootstrapping are inconsistent. In this section, we show that the
corresponding $m$ out of $n$ bootstrap procedures are weakly consistent.
\begin{theorem}\label{thm:m_n_boot1} If $\hat F_n = \mathbb{F}_n$,
$\hat f_n = \tilde f_n$, and $m_n = o(n)$ then the
bootstrap procedure is weakly consistent, for example, (\ref
{eq:boot_limit}) holds in probability.
\end{theorem}
\begin{pf}
Conditions (\ref{eq:cndtn0}), (\ref{eq:cndtn1}) and (\ref{eq:cndtn4})
hold a.s. from (\ref{eq:lil}), as
explained in Section \ref{remarks}. To verify (\ref{eq:cndtn3}), let
$\gamma> 0$ be given. From the proof of Proposition \ref{prop:loc}
[also see Kim and Pollard (\citeyear{kimPo90}), page 218], there exists
$\delta>0$
such that $|\mathbb{F}_{n}(t_0 + h) - \mathbb{F}_{n}(t_0) - F(t_0 + h)
- F(t_0)| \le\gamma h^2 + \mathcal{C}_n n^{-2/3}$,
for $|h| \le\delta$, where $\mathcal{C}_n$'s are random variables of
order $O_P(1)$. We can also assume that $|F(t_0 + h) + F(t_0) - f(t_0)h
- (1/2) f'(t_0) h^2| \le(1/2) \gamma h^2$ for $|h| \le\delta$. Then,
using the inequality $2|a b| \le\gamma a^2 + b^2/\gamma$,
%
%e4.6 ###
%
\begin{eqnarray}\label{eq:A_nB_nC_n}\qquad
&& \biggl| \mathbb{F}_{n}(t_0 + h) - \mathbb{F}_{n}(t_0)
- h \tilde f_n(t_0) - \frac{1}{2}
h^2 f'(t_0) \biggr| \nonumber\\
&&\qquad \le \biggl| \mathbb{F}_{n}(t_0 + h) - \mathbb
{F}_{n}(t_0) - h f(t_0) - \frac{1}{2}
h^2 f'(t_0) \biggr| + |h| |\tilde f_n(t_0) - f(t_0) |
\nonumber\\[-8pt]\\[-8pt]
&&\qquad \le \biggl\{\gamma h^2 + \mathcal{C}_n n^{-{2 /
3}} + \frac{1}{2} \gamma h^2 \biggr\} + \biggl\{ \frac{1}{2} \gamma h^2
+ \frac{1}{2 \gamma} |\tilde
f_n(t_0) - f(t_0) |^2 \biggr\}\nonumber\\
&&\qquad \le 2 \gamma h^2 + \mathcal{C}_n n^{-{2 /3}} +
O_P(n^{-2/3}) \le2 \gamma h^2 + o_P(m_n^{-{2 /3}}).\nonumber
\end{eqnarray}
For (\ref{eq:cndtn2}), write
%
%e4.7 ###
%
\begin{eqnarray}
\label{eq:m_n_boot_t2}
&& m_n^{{2/3}} \{\mathbb{F}_n(t_0 + m_n^{-{1/3}}h) - \mathbb
{F}_n(t_0) - m_n^{-{1/3}}\tilde
f_n(t_0)h\} \nonumber\\
&&\qquad = m_n^{{2/3}} \{ (\mathbb{F}_n - F)(t_0 + m_n^{-{1/3}}h) -
(\mathbb{F}_n - F)(t_0) \}
\nonumber\\[-8pt]\\[-8pt]
&&\qquad\quad{} + m_n^{{1/3}} [f(t_0)-\tilde{f}_n(t_0)]h + \tfrac{1}{
2}f'(t_0)h^2 + o(1) \nonumber\\
&&\qquad \stackrel{P}{\rightarrow} \tfrac{1}{2}f'(t_0)h^2\nonumber
\end{eqnarray}
uniformly on compacts using Hungarian Embedding to bound the second
line and (\ref{eq:chrnff}) (and a two-term Taylor
expansion) in the third.

Given any subsequence $\{n_k\} \subset\mathbb{N}$, there exists a
further subsequence $\{n_{k_l}\}$ such that
(\ref{eq:A_nB_nC_n}) and (\ref{eq:m_n_boot_t2}) hold a.s. and
Theorem \ref{thm:cons_sm_boot} is applicable. Thus,
(\ref{eq:boot_limit}) holds for the subsequence $\{n_{k_l}\}$, thereby
showing that (\ref{eq:boot_limit}) holds in
probability.
\end{pf}

Next consider bootstrapping from $\tilde{F}_n$. We will assume slightly
stronger conditions on $F$, namely, conditions (a)--(d) mentioned in
Theorem 7.2.3 of Robertson, Wright and Dykstra (\citeyear{RWD88}):
\begin{itemize}
\item[(a)] $\alpha_1(F) = \inf\{x\dvtx F(x) = 1\} < \infty$,

\item[(b)] $F$ is twice continuously differentiable on $(0,\alpha_1(F))$,

\item[(c)] $\gamma(F) = \frac{{\sup_{0 < x < \alpha_1(F)}} |f'(x)|}{\inf
_{0 < x < \alpha_1(F)} f^2(x)} < \infty$,

\item[(d)] $\beta(F) = {\inf_{0 < x < \alpha_1(F)}} |\frac
{-f'(x)}{f^2(x)}| > 0$.
\end{itemize}
\begin{theorem}\label{thm:m_n_boot2} Suppose that \textup{(a)--(d)} hold. If
$\hat{F}_n = \tilde{F}_n$, $\hat f_n = \tilde f_n$, and $m_n = o[n
(\log n)^{-{3/2}}]$ then (\ref{eq:boot_limit}) holds in probability.
\end{theorem}
\begin{pf} Conditions (\ref{eq:cndtn0}), (\ref{eq:cndtn1}) and (\ref
{eq:cndtn4}) again follow from (\ref{eq:lil}), as explained in
Section \ref{remarks}. The verification of (\ref{eq:cndtn3}) is
similar to the argument in the proof of Theorem \ref{thm:m_n_boot1}. We
show that (\ref{eq:cndtn2}) holds. Adding and subtracting $m_n^{{2/
3}} [ \mathbb{F}_n(t_0 + m_n^{-{1/3}} h) - \mathbb{F}_n(t_0)]$
from $\mathbb{Z}_{n,2}(h)$ and using
(\ref{eq:m_n_boot_t2}) and the result of Kiefer and Wolfowitz (\citeyear{KW76})
\begin{eqnarray*}
\sup_{|h| \le c} \biggl| {\mathbb Z}_{n,2}(h) - {1\over2}f'(t_0)h^2
\biggr| & \le& 2 m_n^{2/3} \| \tilde{F}_n - {\mathbb F}_n\|
+ o_P(1) \\
& \le& 2 m_n^{2 /3} \| \tilde{F}_n - {\mathbb F}_n\| + o_P(1)
\\
& = & O_P[m_n^{2/3} n^{-{2/3}}\log(n)] + o_P(1)
\end{eqnarray*}
for any $c > 0$ from which (\ref{eq:cndtn2}) follows easily.
\end{pf}

%s5 ###
\section{Discussion}\label{discussion}

We have shown that bootstrap estimators are inconsistent when bootstrap
samples are drawn from either the EDF ${\mathbb F}_n$ or its least
concave majorant ${\tilde F}_n$ but consistent when the bootstrap
samples are drawn from a smoothed version of $\tilde{F}_n$ or an $m$
out of $n$ bootstrap is used. We have also derived necessary conditions
for the bootstrap estimator to have a conditional weak limit, when
bootstrapping from either ${\mathbb F}_n$ or ${\tilde F}_n$ and
presented compelling numerical evidence that these conditions are not
satisfied. While these results have been obtained for the Grenander
estimator, our results and findings have broader implications for the
(in)-consistency of the bootstrap methods in problems with an $n^{1/3}$
convergence rate.

To illustrate the broader implications, we contrast our finding with
those of Abrevaya and Huang (\citeyear{AH05}), who considered a more
general framework, as in Kim and Pollard (\citeyear{kimPo90}). For
simplicity, we use the same notation as in Abrevaya and Huang
(\citeyear{AH05}). Let $W_n := r_n (\theta _n - \theta_0)$ and $\hat
W_n := r_n (\hat\theta_n - \theta_n)$ be the sample and bootstrap
statistics of interest. In our case $r_n = n^{1/3}$, $\theta_0 =
f(t_0)$, $\theta_n = \tilde f_n(t_0)$ and $\hat \theta_n = \tilde
f_n^*(t_0)$. When specialized to the Grenander estimator, Theorem $2$
of Abrevaya and Huang (\citeyear{AH05}) would imply [by calculations
similar to those in their Theorem~5 for the NPMLE in a binary choice
model] that
\[
\hat W_n \Rightarrow\arg\max\hat Z(t) - \arg\max Z(t)
\]
conditional on the original sample, in $P^\infty$-probability, where
$Z(t) = W(t) - ct^2$ and $\hat Z(t) = W(t) + \hat W(t) - ct^2$, $W$ and
$\hat W$ are two independent two sided Brownian motions on $\mathbb{R}$
with $W(0) = \hat W(0) = 0$ and $c$ is a positive constant depending on
$F$. We also know that $W_n \Rightarrow\arg\max Z(t)$
unconditionally. By (v) of Theorem \ref{thm:bootmle}, this would
force the independence of $\arg\max Z(t)$ and $\arg\max\hat Z(t) -
\arg\max Z(t)$; but, there is overwhelming numerical evidence that
these random variables are correlated.

\begin{appendix}\label{app}
%s6 ###
\section*{Appendix}

\begin{lemma}\label{lemma:slope_lcm} Let $\Psi\dvtx\mathbb{R}
\rightarrow
\mathbb{R}$ be a function such that $\Psi(h) \le M$ for all $h \in
\mathbb{R}$, for some $M > 0$, and
%
%e6.1 ###
%
\begin{equation}
\label{eq:Psi_prop} \lim_{|h| \rightarrow
\infty} \frac{\Psi(h)}{|h|} = -\infty.
\end{equation}
Then for any $b > 0$, there exists $c_0 > b$ such that for any $c \ge
c_0$, $L_{\mathbb{R}} \Psi(h) = L_{[-c,c]} \Psi(h)$ for all $|h| \le b$.
\end{lemma}
\begin{pf} Note that for any $c > 0$, $L_{\mathbb{R}} \Psi(h)
\ge L_{[-c,c]} \Psi(h)$ for all $h \in[-c,c]$. Given $b > 0$,
consider $c > b$ and $\Phi_c(h) = L_{[-c,c]} \Psi(h)$ for $h \in
[-b,b]$, and let $\Phi_c$ be the linear extension of $L_{[-c,c]} \Psi
|_{[-b,b]}$ outside $[-b,b]$. We will show that there
exists $c_0 > b + 1$ such that $\Phi_{c_0} \ge\Psi$. Then $\Phi_{c_0}$
will be a concave function everywhere greater than $\Psi$, and thus
$\Phi_{c_0} \ge L_{\mathbb{R}} \Psi$. Hence, $L_{\mathbb{R}} \Psi(h)
\le\Phi_{c_0}(h) = L_{[-c_0,c_0]} \Psi(h)$ for $h \in[-b,b]$,
yielding the desired result.

For any $c > b + 1$, $\Phi_c(h) = \Phi_c(b) - \Phi_c'(b) + \Phi_c'(b)
(h - b + 1)$ for $h \ge b$. Using the min--max formula,
\begin{eqnarray*}
\Phi_c'(b) & = & \min_{-c \le s \le b} \max_{b \le t \le c}
\frac{\Psi(t) - \Psi(s)}{t - s} \\
& \ge& \min_{-c \le s \le b} \frac{\Psi( b + 1) - \Psi(s)}{ (b + 1) -
s} \\
&\ge& \Psi(b + 1) - M =: B_0 \le0.
\end{eqnarray*}
Thus,
\begin{eqnarray*}
\Phi_c(h) & = & \Phi_c(b) - \Phi_c'(b) + \Phi_c'(b) (h - b + 1)
\\
& \ge& \{\Psi(b) - \Phi_c'(b)\} +
\Phi_c'(b) (h - b + 1) \\
& \ge& \Psi(b) + (h - b) B_0
\end{eqnarray*}
for $h \ge b + 1$. Observe that $B_0$ does not depend on $c$. Combining
this with a similar calculation for $h < -(b + 1)$, there are $K_0 \ge
0$ and $K_1 \ge0$, depending only on $b$, for which $\Phi_c(h) \ge K_0
- K_1 |h|$ for $|h| \ge b + 1$. From (\ref{eq:Psi_prop}), there is $c_0
> b + 1$ for which $\Psi(h) \le K_0 - K_1 |h|$ for all $|h| \ge c_0$ in
which case $\Psi(h) \le\Phi_{c_0} (h)$ for all $h$. It follows that
$L_{\mathbb R} \Psi\le\Phi_{c_0}(h)$ for $|h| \le b$.
\end{pf}
\begin{lemma}\label{lem:abc}
Let ${\mathbb B}$ be a standard Brownian motion. If $a, b, c > 0,
a^3b = 1$, then
%
%e6.2 ###
%
\begin{equation}\label{eq:abc1}
P \biggl[\sup_{t \in\mathbb{R}} {|{\mathbb B}(t)|\over a+bt^2} > c
\biggr] = P\biggl[ \sup_{s \in\mathbb{R}} {|{\mathbb B}(s)|\over1+s^2}
> c \biggr].
\end{equation}
\end{lemma}
\begin{pf} This follows directly from rescaling properties of Brownian
motion by letting $t = a^2s$.
\end{pf}
\begin{pf*}{Proof of Proposition \ref{prop:loc}}
Let $J = [a_1,a_2]$ and $\varepsilon> 0$ be as in the statement of the
proposition; let $\gamma= |f'(t_0)|/16$; and recall (\ref
{eq:kmt}) and (\ref{eq:boots_proc1}) from the proof of Proposition \ref
{prop:Z_conv}. Then there exists $0 < \delta< 1$, $C \ge1$, and $n_0
\ge1$ for which (\ref{eq:cndtn3}) and (\ref{eq:cndtn4}) hold for all
$n \ge n_0$. Let $I_{m_n}^* := [-\delta m_n^{1/3}, \delta
m_n^{1/3}]$. By making $\delta$ smaller,\vspace*{-1pt} if necessary, and using
Lemma \ref{lem:loc}, $L_{I_{m_n}}{\mathbb Z}_n(h) =
L_{I_{m_n}^*}{\mathbb Z}_n(h)$ for $|h| \le\delta m_n^{1/3}/2 $
for all but a finite number of $n$ w.p. 1. By increasing the values of
$C$ and $n_0$, if necessary, we may suppose that the right-hand side of
(\ref
{eq:abc1}) (with $c=C$) is less than $\varepsilon/3$, that $P[|\eta| > C]
+ P[\sup_{0 \le t \le1} m_n^{1/6}|{\mathbb E}_{m_n}(t)-{\mathbb
B}_{m_n}^0(t)| > C] \le\varepsilon/3$, and that $L_{I_{m_n}}{\mathbb Z}_n
= L_{I_{m_n}^*}{\mathbb Z}_n$ on $[-{1\over2}\delta m_n^{1/3},
{1\over2}\delta m_n^{1/3}]$ with probability at least $1 -
\varepsilon/3$ for all $n \ge n_0$. We can also assume that $\alpha:=
8C^3/\gamma> 1$. Then, using Lemma \ref{lem:abc} with $a = \alpha
m_n^{-{1/6}}$ and $b = a^{-3}$, the following relations hold
simultaneously with probability at least $1-\varepsilon$ for $n \ge n_0$:
\begin{eqnarray*}
|{\mathbb B}_{m_n}[F_n(t_0) + s] - {\mathbb B}_{m_n}[F_n(t_0)]|
& \le &
C\bigl(\alpha m_n^{-{1/6}} + \alpha^{-3}\sqrt{m_n}s^2\bigr)\qquad \mbox{for
all } s, \\
L_{I_{m_n}}{\mathbb Z}_n & = & L_{I_{m_n}^*}{\mathbb Z}_n \qquad\mbox{on }
\biggl[-\frac{\delta}{2} m_n^{1/3},\frac{\delta}{2} m_n^{1/3}\biggr],
|\eta| \le C,
\end{eqnarray*}
and
\[
\sup_{0 \le t \le1} m_n^{1/6}|
{\mathbb E}_{m_n}(t) - {\mathbb B}_{m_n}^0(t)| \le C.
\]
Let $B_n$ be the event that these four conditions hold. Then $P(B_n)
\ge1-\varepsilon$ for $n \ge n_0$, and from (\ref{eq:boots_proc1}),
$B_n$ implies
%
%e6.3 ###
%
\begin{eqnarray}\label{eq:boundZn1Crude}
|{\mathbb Z}_{n,1}(h)| & \le & C \{\alpha+ \alpha^{-3}m_n^{2/3}
[F_n(t_0+m_n^{-{1/3}}h) - F_n(t_0)]^2 \} + 2C \nonumber\\
&&{} + Cm_n^{1/6}|F_n(t_0 + m_n^{-{1/3}}h) -
F_n(t_0)| \\
&\le& 4 C \{\alpha+ \alpha^{-1}m_n^{2/3}[F_n(t_0+m_n^{-{1/3}}h) -
F_n(t_0)]^2 \}\nonumber
\end{eqnarray}
using the inequalities $|F_n(t_0 + m_n^{-{1/3}} h) - F_n(t_0)| \le
\alpha m_n^{-{1/6}} + \alpha^{-1}m_n^{1/6}[F_n(t_0 +
m_n^{-{1/3}} h) - F_n(t_0)]^2$ and $\alpha> 1$. For sufficiently
large $n$, using (\ref{eq:cndtn4}), we have
%
%e6.4 ###
%
\begin{eqnarray}\label{eq:boundZn1}
|\mathbb{Z}_{n,1} (h)| & \le& 4 C [\alpha+ \alpha^{-1}C^2
m_n^{2/3} (m_n^{-{1/3}}|h| + m_n^{-{1/3}})^2]
\nonumber\\
& \le& 4 C [\alpha+ 2\alpha^{-1} C^2 (h^2 + 1) ] \\
& = & \gamma h^2 + \mathcal{C}\nonumber
\end{eqnarray}
for $|h| \le\delta m_n^{1/3}$ with $\mathcal{C} = 4C\alpha+
8C^3\alpha^{-1}$. Also, we can show that $|{\mathbb Z}_{n,2}(h) -
f'(t_0)h^2 /2| \le\gamma h^2 + \mathcal{C}$
for all $|h| \le\delta m_n^{1/3}$ by (\ref{eq:cndtn3}). Let $b_2
> a_2$ be such that $- 5 \gamma(a_2 + b_2)^2 + 6 \gamma(a_2^{2} +
b_2^2) - 8\mathcal{C} > 0$.

Recalling that $\gamma= -f'(t_0)/16$, $B_n$ implies
\[
-10 \gamma h^2 - 2\mathcal{C} \le{\mathbb Z}_n(h) = {\mathbb
Z}_{n,1}(h) + {\mathbb Z}_{n,2}(h) \le- 6 \gamma h^2 + 2\mathcal{C}
\]
for $|h| \le\delta m_n^{1/3}$ and sufficiently large $n$. Since
the right-hand side is concave, $B_n$ also implies
$L_{I_{m_n}^*}{\mathbb
Z}_n(h) \le- 6\gamma h^2 + 2\mathcal{C}$ for $|h| \le\delta m_n^{1/3}$. Therefore, for
sufficiently large $n$, using the upper bound on $L_{I_{m_n}^*}{\mathbb
Z}_n$, the lower bound on $\mathbb{Z}_n$ obtained above, and
$L_{I_{m_n}}{\mathbb Z}_n(h) = L_{I_{m_n}^*}{\mathbb Z}_n(h)$ for $|h|
\le\delta m_n^{1/3}/2$ on $B_n$, and $[a_2,b_2] \subset I_{m_n}^*$, we have
\begin{eqnarray*}
&&2{\mathbb Z}_n \biggl( {a_2 + b_2 \over2} \biggr) -
[L_{I_{m_n}}{\mathbb Z}_n(a_2) + L_{I_{m_n}}{\mathbb
Z}_n(b_2) ] \\
&&\qquad\ge- 5 \gamma(a_2 + b_2)^2 + 6 \gamma(a_2^{2} + b_2^2) - 8\mathcal
{C} > 0
\end{eqnarray*}
with probability at least $1-\varepsilon$. Thus, $B_n$ implies $2
{\mathbb
Z}_n[{1\over2}(a_2 + b_2)] > L_{I_{m_n}}{\mathbb Z}_n(a_2) +
L_{I_{m_n}}{\mathbb Z}_n(b_2)$ with probability at least $1-\varepsilon$.
Similarly, $B_n$ implies that there is a $b_1 < a_1$ for which $2
{\mathbb Z}_n[{1\over2}(a_1 + b_1)] > L_{I_{m_n}}{\mathbb Z}_n(a_1) +
L_{I_{m_n}}{\mathbb Z}_n(b_1)$ with probability at least $1-\varepsilon$.
Relation (\ref{eq:loc}) then follows from Lemma \ref{lem:ww}. It is
worth noting as a \textit{remark} that $b_1, b_2$ do not depend on the
sequence $F_n$.

Next, consider (\ref{eq:loc2}). Given a compact $J = [-b,b]$, let
$c_{0}(\omega)$ be the smallest positive integer such that for any $c
\ge c_{0}$, $L_{\mathbb{R}} \mathbb{Z}(h) = L_{[-c,c]} \mathbb{Z} (h) $
for $h \in J$. That $c_{0}$ exists and is finite w.p. 1 follows from
Lemma \ref{lemma:slope_lcm}. Defining $W_c := L_{[-c,c]} \mathbb{Z}$
and $Y = L_{\mathbb{R}} \mathbb{Z}$, the event $\{W_c \ne Y \mbox{ on }
J\} \subset\{c_{o} > c\}$. Now given any $\varepsilon> 0$, there exist
$c$ such that $P[c_{o} \le c] > 1 - \varepsilon$. Therefore,
\[
P\bigl[L_{\mathbb{R}} \mathbb{Z} = L_{[-c,c]} \mathbb{Z} \mbox{ on } J \bigr] \ge
P[c_{o} \le c] > 1 - \varepsilon.
\]
\upqed\end{pf*}
\begin{pf*}{Proof of Proposition \ref{prop:edflcm}} First, consider
${\mathbb F}_n$. Let $0 < \gamma< |f'(t_0)|/2$ be given. There is a $0
< \delta< {1\over2}t_0$ such that
%
%e6.5 ###
%
\begin{eqnarray}
\label{eq:2diffF}
\bigl| F(t_0 + h) - F(t_0) - f(t_0)h - \tfrac{1}{2} f'(t_0) h^2 \bigr|
\le\tfrac{1}{2} \gamma h^2
\end{eqnarray}
for $|h| \le2\delta$. From the proof of Proposition \ref{prop:loc},
using arguments similar to deriving (\ref{eq:boundZn1Crude}) and (\ref
{eq:boundZn1}), we can show that
\[
|({\mathbb F}_n-F)(t_0+h)-({\mathbb F}_n-F)(t_0)| < \tfrac{1}{2}\gamma
h^2 + C n^{-{2 /3}}
\]
for $|h| \le2 \delta$ with probability at least $1 - \varepsilon$ for
sufficiently large $n$. Therefore, by adding and subtracting $F(t_0+h)
- F(t_0)$ and using (\ref{eq:2diffF}),
%
%e6.6 ###
%
\begin{equation}
\label{eq:edf}
\bigl| {\mathbb F}_n(t_0 + h) - {\mathbb F}_n(t_0) - f(t_0)h - \tfrac
{1}{2} f'(t_0) h^2 \bigr|
\le\gamma h^2 + C n^{-{2/3}}
\end{equation}
for $|h| \le2 \delta$ with probability at least $1-\varepsilon$ for
large $n$.

Next, consider $\tilde{F}_n$. Let $B_n$ denote the event that (\ref
{eq:edf}) holds. Then $P(B_n)$ is eventually
larger than $1 - \varepsilon$ and on $B_n$, we have
\[
{\mathbb F}_n(t_0+h) - \mathbb{F}_n(t_0) - f(t_0) h \le\bigl\{ \gamma
- \tfrac{1}{2} |f'(t_0)| \bigr\} h^2 + Cn^{-{2/3}}
\]
for $|h| \le2 \delta$. Let $E_n$ be the event that $\tilde{F}_n(h) =
L_{[t_0-2\delta,t_0+2\delta]}{\mathbb F}_n(h)$ for $h \in[t_0-
\delta,t_0+ \delta]$. Then by Lemma \ref {lem:loc}, $P(E_n) \ge1 -
\varepsilon$, for all sufficiently large $n$. Taking concave majorants
on either side of the above display for $|h| \le2 \delta$ and noting
that the right-hand side of the display is already concave, we have:
${\tilde F}_n(t_0+h) - \mathbb{F}_n(t_0) - f(t_0) h \le\{ \gamma-
\frac{1}{2} |f'(t_0)| \} h^2 + C n^{-{2/3}}$, for $|h| \le\delta$ on
$B_n \cap E_n$. Setting $h = 0$ shows that on $E_n \cap B_n$, $\tilde
F_n(t_0) - \mathbb{F}_n(t_0) \le C n^{-{2/3}}$. Now, as ${\mathbb
F}_n(t_0) \le\tilde{F}_n(t_0)$, it is also the case that on $E_n \cap
B_n$, for $|h| \le\delta$,
%
%e6.7 ###
%
\begin{equation}
\label{eq:cncvbnd2}
{\tilde F}_n(t_0+h) - \tilde F_n(t_0) - f(t_0) h \le\bigl\{ \gamma-
\tfrac{1}{2} |f'(t_0)| \bigr\} h^2 + C n^{-{2/3}}.
\end{equation}
Furthermore on $E_n \cap B_n$,
%
%e6.8 ###
%
\begin{eqnarray}
\label{eq:cncvbnd3}
&& \tilde F_n(t_0 + h) - \tilde F_n(t_0) - f(t_0) h - \tfrac{1}{2}
f'(t_0) h^2 \nonumber\\
&&\qquad \ge \mathbb F_n(t_0 + h) - \{\mathbb F_n(t_0) + C n^{-{2/3}}\}
- f(t_0) h - \tfrac{1}{2} f'(t_0) h^2
\\
&&\qquad \ge -\gamma h^2 - 2 C n^{-{2/3}}.\nonumber
\end{eqnarray}
Therefore, combining (\ref{eq:cncvbnd2}) and (\ref{eq:cncvbnd3}),
\[
\bigl| {\tilde F}_n(t_0 + h) - {\tilde F}_n(t_0) - f(t_0)h - \tfrac
{1}{2} f'(t_0) h^2 \bigr| \le\gamma h^2 + 2C n^{-{2/3}}
\]
for $|h| \le\delta$ with probability at least $1 - 2\varepsilon$ for
large $n$.
\end{pf*}
\end{appendix}

\printaddresses

\end{document}